\documentclass[12pt]{article}		

\usepackage{amsmath}
\usepackage{graphicx}
\usepackage{amsxtra}
\usepackage{amstext}
\usepackage{amssymb}
\usepackage{amsthm}
\usepackage{latexsym}
\newtheorem{theorem}{Theorem}[section]
\newtheorem{lemma}{Lemma}[section]
\newtheorem{definition}{Definition}[section]

\alph{footnote}				
\title{Prewavelet Solution to Poisson Equations}

\author{Ming-Jun Lai\footnote{mjlai@math.uga.edu. This author is partly
supported by the National Science Foundation under grant EAR-0327577 } 
 and Haipeng Liu\\
Department of Mathematics\\ University of Georgia\\ Athens, GA 30602}

\date{\today}				

\begin{document}			

\pagenumbering{arabic}			
\maketitle				
\thispagestyle{empty}			

\begin{abstract}	
Finite element method is one of powerful numerical methods to solve PDE. Usually, 
if a finite element solution to a Poisson equation based on a triangulation 
of the underlying domain is not accurate enough,  one will discard the solution and 
then refine 
the triangulation uniformly and  compute a new finite element solution over 
the refined triangulation.  
It is wasteful to discard the original finite element solution. 
We propose a prewavelet method to save  the original solution by adding a 
prewavelet subsolution to obtain  the refined level finite element solution. 
To increase the accuracy of numerical solution to Poisson equations, we can keep adding 
prewavelet subsolutions. 

Our prewavelets are orthogonal in the $H^1$ norm and 
they are compactly supported except for one globally 
supported basis function in a rectangular domain.   
We have implemented these prewavelet basis functions in MATLAB and used them for numerical 
solution of Poisson equation with Dirichlet boundary conditions. 
Numerical simulation demonstrates that our prewavelet solution is much more 
efficient than the standard finite element method. 

\end{abstract}	

\section{ Introduction }	
Finite element method is one of powerful numerical methods to solve PDE. Usually, 
if a finite element solution to a Poisson equation based on one level triangulation 
of the underlying domain is not accurate enough,  one will discard the solution and then 
refine the triangulation and  compute a new finite element solution at the refined level. 
It is wasteful to throw the original finite element solution away. 
In order to save the original
solution and get the more accurate new solution, we have to add $H^1$ orthogonal subsolution.
That is, let $V_h$ be a finite element space over a triangulation $\Delta_h$ 
and $V_{h/2}$ be the finite element space over the refined triangulation. Since $V_h\subset
V_{h/2}$, let $W_h= V_{h/2}\ominus V_h$ under $H^1$ norm, if $\Phi_h\in V_h$
is a finite element solution of Poisson equation with Dirichlet boundary condition, we can
find $\Psi_h\in W_h$ so that $\Phi_h+\Psi_h$ is the finite element solution in $V_{h/2}$.
In addition, suppose that 
$\phi_h$ is the most accurate solution that a computer can compute in the sense that it would
be out of memory when computing a finite element solution $\Phi_{h/2}$ in $V_{h/2}$ directly. 
Since the size of the linear system associated with $\Psi_h$ is smaller than $\Phi_{h/2}$,
if the computer can solve $\Psi_h$, we can add $\Psi_h$ to $\Phi_h$ to get $\Phi_{h/2}$ 
achieving the next level of accuracy.   
In this paper, we discuss how to compute $\Psi_h$.  We shall construct  compactly
supported basis functions and a global supported basis function $\psi_{h,k}, k=1, \cdots, N_h$
which span $W_h$. $\psi_{h,k}$'s are called prewavelets and $\Psi_h$ is a linear combination
of these $\psi_{h,k}$'s and hence is called a prewavelet subsolution.

Prewavelets have been studied for more than 10 years (cf. \cite{JM90}, \cite{CSW92}). 
There are many methods available 
to construct compactly supported prewavelets over 2D domains 
under the $L_2$ norm. That is $W_h= V_{h/2}\ominus V_{h}$ under $L_2$ norm, e.g., in a
series of papers \cite{FQ99}, \cite{FQ01}, \cite{HM99}, \cite{KO95}, and \cite{BDG01}. 
In 1997, Bastin and Laubin (\cite{BL97}) explained how to construct 
compactly supported orthonormal wavelets in Sobolev space in the univariate setting. 
See also \cite{BB98} for biorthogonal wavelets in Sobolev space. 
In \cite{LO97}, Lorentz and Oswald showed that there is 
no compactly supported prewavelets in Sobolev space or under $H^1$ norm based on 
integer translations of a box spline over ${\bf R}^2$. Since continuous piecewise linear 
finite element can be expressed by using box spline $B_{111}$, the result in \cite{LO97}
ruins a hope to find compactly supported prewavelets under $H^1$ norm. 
But this is not an end of story. It is possible to construct compactly supported 
prewavelets in a semi-norm in the univariate setting in \cite{JWZ03}. It is also possible to
construct compactly supported prewavelets in $H^r$ norm over each nested subspace, but
the union of these prewavelets over all levels fails to be a stable basis for a Sobolev 
space (cf. \cite{L06}).  Our new question is if we can find a prewavelet basis 
with as few as possible global supported prewavelet functions.  
Our anwser is affirmative. That is, there is a prewavelet basis for $W_h$ with only 
one global supported basis function under the $H^1$ norm over rectangular domains. Also
it is possible to find a compactly supported prewavelet basis for $W_h$ under the $H^1$ norm 
for Poisson equation over a triangular domain (cf. \cite{Liu07}). 

The paper is organized as follows: We first explain that the Dirichlet boundary value problem
of Poisson equation can be converted into a Poisson equation with zero boundary condition. 
An explicit conversion will be given. 
Thus the $H^1$ norm is now equivalent to the $H^1_0$ semi-norm. Then we introduce some notation
to explain the weak solution of Poisson equation and its approximation to
the exact solution.  These explanations are well-known and 
given in the Preliminary section \S 2. In \S 3, we explain how to construct compactly 
supported prewavelets under $H^1_0$ semi-norm. In \S 4, we explain how to implement our 
prewavelet method for numerical solution of Poisson equation. Finally in \S 5 we present 
some numerical results. Our numerical experiment show that the time for computing a finite 
element solution by our prewavelet method is about half of the time by 
the standard finite element method using the direct method for inverting the linear systems.
If using the conjugate gradient method for the linear systems for the finite element method,
the prewavelet method is still faster than for sufficiently accurate iterative solutions.

\section{Preliminary} 
Let us start with a square domain $\Omega=(0,1)\times(0,1)\in R^2$. 
Consider the Dirichlet boundary value problem for Poisson equation:
\begin{equation*}
\begin{split}
     \left\{ 
          \begin{array}{ll}
  \quad-\Delta u(x,y)=g(x,y) , \hspace{0.5cm} (x,y)\in \Omega\\
  \quad u(x,y)=f_1(x), \quad for \quad y=0 \quad and \quad 0\leq x\leq 1\\
  \quad u(x,y)=f_2(x), \quad for \quad y=1 \quad and \quad 0\leq x\leq 1\\
  \quad u(x,y)=f_3(y), \quad for \quad x=0 \quad and \quad 0\leq y\leq 1\\
  \quad u(x,y)=f_4(y), \quad for \quad x=1 \quad and \quad 0\leq y\leq 1
          \end{array} 
     \right. 
\end{split}  
\end{equation*}
Without lose of generality, 
we may assume that $f_1(1)=f_2(1)=f_3(1)=f_4(1)=f_1(0)=f_2(0)=f_3(0)=f_4(0)=0$.  
Otherwise, letting $f_1(0)=f_3(0)=a_1$, 
$f_3(1)=f_2(0)=a_2$, $f_2(1)=f_4(1)=a_3$, $f_4(0)=f_1(1)=a_4$,  
we define $h(x,y)=a_1+(a_4-a_1)x+(a_2-a_1)y+(a_3+a_1-a_4-a_2)xy$, 
and $v(x,y)=u(x,y)-h(x,y)$. Then the above Dirichlet problem becomes to:
\begin{equation*}
\begin{split}
     \left\{ 
          \begin{array}{ll}
  \quad-\Delta v(x,y)=g(x,y) , \hspace{0.5cm} (x,y)\in \Omega\\
  \quad v(x,y)=f_1(x)-h(x,0) , \quad for \quad y=0 \quad and \quad 0\leq x\leq 1\\
  \quad v(x,y)=f_2(x)-h(x,1) , \quad for \quad y=1 \quad and \quad 0\leq x\leq 1\\
  \quad v(x,y)=f_3(y)-h(0,y) , \quad for \quad x=0 \quad and \quad 0\leq y\leq 1\\
  \quad v(x,y)=f_4(y)-h(1,y) , \quad for \quad x=1 \quad and \quad 0\leq y\leq 1
          \end{array} 
     \right. 
\end{split}  
\end{equation*}
which satisfy the above assumption.

Now let $w(x)=v(x,y)-x(f_4(y)-h(1,y))-(1-x)(f_3(y)-h(0,y))
-y(f_2(x)-h(x,1))-(1-y)(f_1(x)-h(x,0))$. Then $w(x)$ satisfies the equation 
\begin{equation*}
\begin{split}
     \left\{ 
          \begin{array}{ll}
  \quad-\Delta w(x,y)=g_1(x,y) , \hspace{0.5cm} & (x,y)\in \Omega\\
  \quad w(x,y)=0, \hspace{0.5cm} & (x,y)\in \partial \Omega \quad 
          \end{array} 
     \right. 
\end{split}  
\end{equation*}
with \( g_1(x,y)= g(x,y)+\frac{\partial^2}{\partial y^2}
\lbrack-x(f_4(y)-h(1,y))-(1-x)(f_3(y)-h(0,y))\rbrack
        +\frac{\partial^2}{\partial x^2} \lbrack-y(f_2(x)-h(x,1))-(1-y)
(f_1(x)-h(x,0))\rbrack\).

If we can find solution for $w$, it is easy to get $u(x,y)$.
In the remaining paper, 
we only consider the Poisson equation with zero boundary condition:
\begin{equation}
     \left\{ 
          \begin{array}{ll}
  \quad-\Delta u(x,y)=g(x,y) , \hspace{0.5cm} & (x,y)\in \Omega\\
  \quad u(x,y)=0, \hspace{0.5cm} & (x,y)\in \partial \Omega .
          \end{array} 
     \right. 
\label{PDE}
\end{equation}
Next we define 
$$
H^1_0(\Omega)=\{v\in L^2(\Omega):\hspace{0.5cm}  
\langle v, v\rangle_s<\infty\hspace{0.5cm}  and \hspace{0.5cm} 
v(x,y)=0, (x,y)\in \partial \Omega\}, 
$$
where  the inner product $\langle u, v\rangle_s$ is defined by
$$ 
\langle u, v \rangle_s=
\int_0^1\int_0^1\frac{\partial u(x,y)}{\partial x}\frac{\partial v(x,y)}{\partial x}
+\frac{\partial u(x,y)}{\partial y}\frac{\partial v(x,y)}{\partial y}dxdy.
$$
By using Poincare's inequality, $\|u\|_s=\sqrt{\langle u, u\rangle_s}$ is 
a standard Sobolev norm for $H^1_0(\Omega)$. Suppose $u, v\in H^1_0(\Omega)$. 
Integration by parts yields
\begin{equation*}
\begin{split}
\langle g, v\rangle &=\int_0^2\int_0^2g(x,y)v(x,y)dxdy\\
                   &=\int_0^2\int_0^2 -\Delta u(x,y)v(x,y)dxdy\\
                   &=\int_0^2\int_0^2 \frac{\partial u(x,y)}{\partial x}
                     \frac{\partial v(x,y)}{\partial x}+
              \frac{\partial u(x,y)}{\partial y}\frac{\partial v(x,y)}{\partial y}dxdy\\
                        &=\langle u, v \rangle_s.
\end{split}  
\end{equation*}
Thus, a weak solution $u$ to (\ref{PDE}) is characterized by finding $u\in H^1_0(\Omega)$ such
that 
\begin{equation}
\langle u,v \rangle_s=\langle g, v\rangle, \quad \forall v\in H^1_0(\Omega).
\label{WPDE}
\end{equation}

The following result  is well-known. For convenience, we present a short proof.  
\begin{theorem}
Suppose $g\in C(\Omega)$ and $u\in C^2(\Omega)$ satisfy (\ref{WPDE}).  
Then u is weak solution of (\ref{PDE}).
\end{theorem}

\begin{proof}
Let $v\in H^1_0(\Omega)$. Then integration by parts gives
\begin{equation*}
\begin{split}
     &\quad \langle g, v\rangle=\langle u, v\rangle_s \\
     &\quad=\int_0^1\int_0^1 \frac{\partial u(x,y)}{\partial x}\frac{\partial v(x,y)}
       {\partial x}+\frac{\partial u(x,y)}{\partial y}\frac{\partial v(x,y)}{\partial y}dxdy\\
     &\quad=\int_0^1\int_0^1 -\Delta u(x,y)v(x,y)dxdy \\
     &\quad=\langle -\Delta u(x,y), v\rangle.
\end{split}
\end{equation*}
It follows that $\langle g-(-\Delta u(x,y)), v\rangle =0$ for all $v\in H^1_0(\Omega)$.  
That is, $g\equiv -\Delta u$ and hence, $u$ satisfies (\ref{PDE}). 
\end{proof}

Next we introduce continuous linear spline space on $\Omega=[0,1]\times [0, 1]$. 
For convenience, let $N_j=(2^j-1)^2$ and $j\ge 1$. Denote 
$x_{ji}=\frac{i}{2^j}=y_{ji}$ for $i=1,..,2^j-1$. Clearly, the lines 
segment of $x=x_{ji}$ and $y=y_{jk}$ divide the square $\Omega$ into $N_j$ sub-squares. 
The diagonal going from down-left to up-right of each sub-square 
divides the sub-square into two congruent triangle. We will refer to the set of 
all such triangles  as a Type-1 triangulation of $\Omega$ (see Figure 1).

\setlength{\unitlength}{.4in}                             

\begin{center}
\begin{picture}(10,5)(0,0)
\linethickness{1pt}
\put(0,1){\line(1,0){4}}
\put(0,2){\line(1,0){4}}
\put(0,3){\line(1,0){4}}
\put(0,4){\line(1,0){4}}
\put(0,5){\line(1,0){4}}

\put(0,1){\line(0,1){4}}
\put(1,1){\line(0,1){4}}
\put(2,1){\line(0,1){4}}
\put(3,1){\line(0,1){4}}
\put(4,1){\line(0,1){4}}

\put(0,4){\line(1,1){1}}
\put(0,3){\line(1,1){2}}
\put(0,2){\line(1,1){3}}
\put(0,1){\line(1,1){4}}
\put(1,1){\line(1,1){3}}
\put(2,1){\line(1,1){2}}
\put(3,1){\line(1,1){1}}
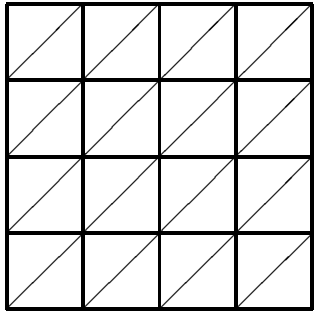
\put(2,0.5){\makebox(0,0){Figure 1. Type-I triangulation with j=2.}}
\end{picture}
\end{center}

Define $\phi^j_{ik}$ to be linear spline with support on the hexagon with following vertices
$$
\{(x_{j(i-1)},y_{j(k-1)}),(x_{ji},y_{j(k-1)}),(x_{j(i+1)},y_{j(k)})
,(x_{j(i+1)},y_{j(k+1)}),(x_{j(i+1)},y_{j(k)}),(x_{j(i-1)},y_{j(k)}) \}
$$
and $\phi_{ik}(x_{ji^{'}},y_{jk^{'}})=\delta_{i,i^{'}}\delta_{k,k^{'}}$, where 
$\delta_{i,i'} =0$ if $i'\not=i$ and $1$ if $i'=i$. 

Let $V_j=span\{\phi^j_{ik},i=1,..,2^j-1,k=1,..,2^j-1\}$ be the subspace 
of $H^1_0(\Omega)$. By following lemma, there exists a unique $u_j\in V_j$ 
satisfying 
\begin{equation}
\langle u_j, v\rangle_s=\langle f,v \rangle \hspace{0.2cm} \forall v \in V_j.
\label{UNI}
\end{equation}
$u_j$ is the standard finite element solution in $V_j$. The following result is well-known.
For completeness, we include a short proof. 

\begin{lemma}
Given $g\in L^2(\Omega)$, (\ref{UNI}) has a unique solution.
\end{lemma}

\begin{proof}
Reorder the basis functions $\phi^{(j)}_{ik}$ to $\phi_m$, $m=1,...,N_j$ and 
let $u_j=\sum a_m\phi_m$. Denote $k_{mn}=\langle \phi_m, \phi_n\rangle_s$ and 
$F_m=\langle f, \phi_m \rangle$ for $m=1,.....,N_j$.
Set $A=(a_m)$ to be the coefficient vector, $K=[k_{mn}]_{1\leq m,n \leq N_j}$ to be the 
stiff matrix, and $F=(F_m)$ to be the right hand side vector. 
Then the solutions in (\ref{UNI}) is written in the following matrix equation form
\begin{equation}
     KA=F.
\label{MEq}
\end{equation}
We claim that the solution for above equation always exists and is unique. Otherwise there is 
a nonzero vector ${\bf c}$ such that $K{\bf c}=0$. 
Write ${\bf c}=(c_m, m=1,......,N_j)$ and let $v=\sum_{i=1}^{N_j}c_i\phi_i$ be the
linear spline. Then $K{\bf c}=0$ is equivalent to 
\begin{equation*}
    \langle v,\phi_m\rangle_s=0 \quad \forall  m=1, \cdots, N_j. 
\end{equation*}
Multiplying $\langle v,\phi_m\rangle_s$ by $c_m$ and 
summing over m yields $\langle v, v\rangle_s=0$. Thus, $v=a + bx+c dxy$. 
Boundary condition implies $v\equiv 0$. 
Since $\{\phi_m\}$ are linear independent, ${\bf c}\equiv0$ and hence, the solution is  
unique.
\end{proof}

Let us discuss the error between $u$ and $u_j$. It is standard in finite element analysis
(cf. \cite{BS94}). For completeness we present a simple derivation. 
Subtracting ($\ref{UNI}$) from 
(\ref{WPDE}) implies
\begin{equation}
\langle u-u_j,w\rangle_s=0 \hspace{0.2cm}\quad\quad\quad \forall w \in V_j .
\end{equation}
Then for any $v\in V_j$
\begin{equation*}
\begin{split}
\|u-u_j\|_s^2 &\quad = \langle u-u_j,u-u_j\rangle_s\\
              &\quad =\langle u-u_j,u-v\rangle_s+\langle u-u_j, v-u_j\rangle_s\\
              &\quad =\langle u-u_j, u-v\rangle_s\\
              &\quad \leq\|u-u_j\|_s\|u-v\|_s\\ 
\end{split}
\end{equation*}
It follows that $\|u-u_j\|_s\leq \|u-v\|_s$ for any $v\in V_j$. 
Thus we have proved the following.

\begin{lemma}(C\'ea's Lemma)
$\|u-u_j\|_s = \min \{\|u-v\|_s : v\in V_j\}.$
\end{lemma}

Given $u\in C^0(\Omega)$, let $u_j\in V_j$ be the interpolant of $v$: 
$$
u_j=\sum_{ik} u(x_{ji},y_{jk})\phi_{ik}^{(j)}.
$$
The following error estimate is well-known.
\begin{lemma}
Suppose $u\in C^2(\Omega)$. Then
$$\|u-u_j\|_s\leq \frac{\sqrt{12}}{2^j}  
\sqrt{\left\|\frac{\partial^2 u}{\partial x^2}\right\|^2_{L^{\infty}}
+\left\|\frac{\partial u}{\partial x}\frac{\partial u}{\partial      
 y}\right\|^2_{L^{\infty}}+\left\|\frac{\partial^2 u}{\partial y^2}\right\|^2_{L^{\infty}}}.$$
\end{lemma}

\begin{proof}
The proof is elementary and is left to the reader. See \cite{Liu07} for detail. 
\end{proof}

\section{Multiresolution and Prewavelets over Type-I triangulations}
We start with the definition of multi-resolution approximation of $H_0^1(\Omega)$: 
\begin{definition}
A multiresolution approximation of $H_0^1(\Omega)$ is a sequence of finite 
dimensions subspaces $V_j$, $j\in Z^{+}$ of $H_0^1(\Omega)$ such that 
\begin{description}
\item{(1)}  $ V_j\subset V_{j+1},\quad j\in Z^{+}$; 
\item{(2)} $ \bigcup_{j=1}^{\infty}V_j$ is dense in $H_0^1(\Omega)$.
\end{description}
\end{definition}

Let $\Gamma^j$ be the type-1 triangulation with $2N_j$ triangles. 
Naturally, let $\Gamma^{j+1}$ be the uniform refinement of $\Gamma^j$. 
Let $V_j$ be the continuous piecewise linear spline space defined on the previous section.
That is,  $V_{j}=span\{\phi^{j}_{ik},i=1,..,2^{j}-1,k=1,..,2^{j}-1\}$, where
$\phi^{j}_{ik}$  are continuous piecewise linear functions which is 1 at $(x_{ji}, 
y_{jk})$ and zero at all other vertices.   
Let $V_{j+1}=span\{\phi^{j+1}_{ik},i=1,..,2^{j+1}-1,k=1,..,2^{j+1}-1\}$, 
and $(x_{j+1,i},y_{j+1,k})$ are the vertices on the j+1 level Type-1 triangulation. 
Then the refinement equation is easily seen to be 
$$\phi^j_{ik}=\phi^{j+1}_{2i,2k}+\frac{1}{2}\phi^{j+1}_{2i-1,2k}+
\frac{1}{2}\phi^{j+1}_{2i-1,2k-1}+\frac{1}{2}\phi^{j+1}_{2i,2k-1}+\frac{1}{2}
\phi^{j+1}_{2i+1,2k}+\frac{1}{2}\phi^{j+1}_{2i+1,2k+1}+\frac{1}{2}\phi^{j+1}_{2i,2k+1}.$$
See the Figure 2.
\begin{center} 
\begin{picture}(6,5)(0,0)
\linethickness{0.5pt}
\put(2,1){\line(1,0){2}}
\put(2,3){\line(1,0){0.8}}
\put(3.2,3){\line(1,0){0.6}}
\put(4.2,3){\line(1,0){0.6}}
\put(5.2,3){\line(1,0){0.8}}
\put(4,5){\line(1,0){2}}

\put(2,1){\line(0,1){2}}
\put(4,1){\line(0,1){0.7}}
\put(4,2.3){\line(0,1){0.4}}
\put(4,3.3){\line(0,1){0.4}}
\put(4,4.3){\line(0,1){0.7}}
\put(6,3){\line(0,1){2}}


\put(2,3){\line(1,1){2}}
\put(2,1){\line(1,1){0.8}}
\put(3.2,2.2){\line(1,1){0.6}}
\put(4.2,3.2){\line(1,1){0.6}}
\put(5.2,4.2){\line(1,1){0.8}}
\put(4,1){\line(1,1){2}}

\put(2.89,2.86){$\frac{1}{2}$}
\put(4.89,2.86){$\frac{1}{2}$}
\put(3.89,3.86){$\frac{1}{2}$}
\put(3.89,1.86){$\frac{1}{2}$}
\put(2.89,1.86){$\frac{1}{2}$}
\put(4.89,3.86){$\frac{1}{2}$}
\put(3.89,2.86){2}

\put(4,0){\makebox(0,0){Figure 2. Dilation relations}}

\put(1.0,0.5){(0,0)}
\put(3.0,0.5){$(1/2^j,0)$}
\put(0.5,3){$(0,1/2^j)$}
\end{picture}
\end{center}

The main purpose of this paper is to build a basis for the orthogonal complement $W_j$ 
of $V_j$ in $V_{j+1}$ under  the inner product $\langle \cdot ,\cdot \rangle_s$.
Suppose we have the $W_j$. Then $V_{j+1}=V_j+W_j$ under the $H^1_0(\Omega)$ inner product.
For a solution $u_j$ satisfying (3), we do not have to find out the solution for 
\begin{equation*}
\begin{split}
u_{j+1}\in V_{j+1} \hbox{ such that } 
\langle u_{j+1},v\rangle_s=\langle g, v \rangle \hspace{0.2cm} \forall v \in V_{j+1}.
\end{split}
\end{equation*}
Instead,  we only need to find solutions for 
\begin{equation*}
\begin{split}
w_{j}\in W_{j} \hbox{ such that }
\langle w_{j},v\rangle_s=\langle g, v \rangle \hspace{0.2cm} \forall v \in W_{j}.
\end{split}
\end{equation*}
Then we have $w_j+u_j=u_{j+1}$. Ideally, we hope the supports of basis functions 
for $W_{j}$ are small, since small support can accelerate the calculations 
of $\langle g, v \rangle_s$. As explained in the Introduction, there is no compactly
supported prewavelets for $W_j$. Neverthless, we shall construct basis functions 
with only one globally supported basis function for $W_j$ in the following. 

Clearly the $\Gamma_{j}$ can be continuously refined and hence
we will have a nested sequence of subspaces 
$$V_{1}\subset V_{2}\subset V_{3}\subset V_{4}\subset V_{5}......$$
to span $H^1_0(\Omega)$ by Lemmma 2.3 since $C^2(\Omega)$ is dense in $H^1_0(\Omega)$.

Let $W_{j}\subset V_{j+1}$ be the orthogonal complement
of $V_j$ in $V_{j+1}$ for  each refinement level $j$, i.e., 
$$V_{j+1}=V_j \bigoplus W_j.$$
Then we get the decomposition
$$ V_{j+1}=V_1 \bigoplus W_1\bigoplus W_2\bigoplus W_3\bigoplus......\bigoplus W_j$$
for any $j\geq 1$. The weak solution $u_{j+1}$ to the Poisson equation (\ref{PDE}) 
at $V_{j+1}$ can be built by 
$$
u_{j+1}=u_1+w_1+w_2+\cdots+ w_j.
$$ 

We now focus on building basis functions for the orthogonal complement $W_j$.
By direct calculation, we obtain the following lemma immediately.
\begin{lemma} We have 
$\langle \phi^{j}_{ik}, \phi^{j+1}_{2i,2k},\rangle_s=2 $,\\
$\langle \phi^{j}_{ik}, \phi^{j+1}_{2i-1,2k},\rangle_s=1/2 $,\quad
$\langle \phi^{j}_{ik}, \phi^{j+1}_{2i,2k-1},\rangle_s=1/2 $,\quad
$\langle \phi^{j}_{ik}, \phi^{j+1}_{2i+1,2k},\rangle_s=1/2 $,\\
$\langle \phi^{j}_{ik}, \phi^{j+1}_{2i,2k+1},\rangle_s=1/2 $,\quad
$\langle \phi^{j}_{ik}, \phi^{j+1}_{2i-1,2k-1},\rangle_s=1 $,\quad
$\langle \phi^{j}_{ik}, \phi^{j+1}_{2i+1,2k+1},\rangle_s=1 $,\\
$\langle \phi^{j}_{ik}, \phi^{j+1}_{2i-2,2k},\rangle_s=-1/2 $,\quad
$\langle \phi^{j}_{ik}, \phi^{j+1}_{2i+2,2k},\rangle_s=-1/2 $,\quad
$\langle \phi^{j}_{ik}, \phi^{j+1}_{2i,2k-2},\rangle_s=-1/2 $,\\
$\langle \phi^{j}_{ik}, \phi^{j+1}_{2i,2k+2},\rangle_s=-1/2 $,\quad
$\langle \phi^{j}_{ik}, \phi^{j+1}_{2i-2,2k-2},\rangle_s=0 $,\quad
$\langle \phi^{j}_{ik}, \phi^{j+1}_{2i+2,2k+2},\rangle_s=0 $,\\
$\langle \phi^{j}_{ik}, \phi^{j+1}_{2i-2,2k-1},\rangle_s=-1/2 $,\quad
$\langle \phi^{j}_{ik}, \phi^{j+1}_{2i-1,2k+1},\rangle_s=-1 $,\quad
$\langle \phi^{j}_{ik}, \phi^{j+1}_{2i+1,2k+2},\rangle_s=-1/2 $,\\
$\langle \phi^{j}_{ik}, \phi^{j+1}_{2i+2,2k+1},\rangle_s=-1/2 $,\quad
$\langle \phi^{j}_{ik}, \phi^{j+1}_{2i+1,2k-1},\rangle_s=-1 $,\quad
$\langle \phi^{j}_{ik}, \phi^{j+1}_{2i-1,2k-2},\rangle_s=-1/2 $,\\
$\langle \phi^{j}_{ik}, \phi^{j+1}_{i',k'},\rangle_s=0 $, 
for other $i^{'},k^{'}$ which are not listed above.
\label{inner}
\end{lemma}

Let $\psi^j$ be a function in $W_j$. Since $W_j\subset V_{j+1}$, let us write 
$\psi^j=\sum_{ik} \phi^{j+1}_{ik}b_{ik}$ for some unknown coefficients $b_{ik}$. 
Then by orthogonal condition $\langle \phi^{j}_{i{'}k{'}}, \psi^j\rangle_s=0 $, 
we need to solve the following equations.

\begin{equation}
0=\langle \phi^{j}_{i{'}k{'}}, \sum_{i,k} b_{ik} \phi^{j+1}_{ik}\rangle_s
=\sum_{i,k}b_{ik} \langle \phi^{j}_{i{'}k{'}}, \phi^{j+1}_{ik}\rangle_s.
\label{WWLT}
\end{equation}
Each $(i^{'}k^{'})$ determines one equation. 
Since there are $N_j$ elements in the set $V_j$, they determine the $N_j$ equations. 
These $N_j$ equations with $N_{j+1}$ coefficients, $b_{i,k}$. There are at least $N_{j+1}
-N_j$ degrees of freedom.  
The solution space of these equation system should be the $W_{j}$. The linear independence
of $\phi^j_{i',k'}$ implies that the coefficient matrix of the above linear system is of
full rank. Hence, 
there are $N_{j+1}-N_{j}$ linear independent solutions 
which constitute a basis for $W_{j}$.

\begin{definition}
Let $V_{j+1}^m=span\{\phi^{j+1}_{ik}, i=1,..,2m-1,k=1,..,2m-1\}$ 
be a subspace of $V_{j+1}$.  Let $W_{j}^m$ be subspace of $W_{j}$ such that 
$W_{j}^m=W_{j}\bigcap V_{j+1}^m$. 
\end{definition}

Obviously $\emptyset\subset V_{j+1}^1\subset V_{j+1}^2\subset \ldots 
\subset V_{j+1}^{2^{j}}=V_{j+1}$, 
and $\emptyset\subset W_{j}^1\subset W_{j}^2\subset \ldots \subset W_{j}^{2^{j}}=W_j$.
There is no nonzero solution of (\ref{WWLT}) in space of $V_{j+1}^1$. However,  
there are five solution of (\ref{WWLT}) in space $V_{j+1}^2$. 
They are solutions of the following system of linear equations.

$$\begin{array}{cc}
\displaystyle \sum_{1\leq i,k\leq 3} b_{ik}\left\langle  \phi^{j+1}_{ik}, 
\phi^j_{1,1}\right\rangle_s=0, & 
\displaystyle \sum_{1\leq i,k\leq 3}b_{ik} \left\langle \phi^{j+1}_{ik}, \phi^j_{2,1}
\right\rangle_s=0, \\
\displaystyle \sum_{1\leq i,k\leq 3}b_{ik} \left\langle \phi^{j+1}_{ik}, \phi^j_{1,2}
\right\rangle_s=0, & 
\displaystyle \sum_{1\leq i,k\leq 3}b_{ik} \left\langle \phi^{j+1}_{ik}, \phi^j_{2,2}\right
\rangle_s=0. 
\end{array}$$
They are equivalent to the following equations.
$$\left(\begin{array}{ccccccc}
\langle \phi^j_{1,1},\phi^{j+1}_{1,1}\rangle_s&&\langle \phi^j_{1,1},\phi^{j+1}_{2,1}\rangle_s&&....&&\langle \phi^j_{1,1},\phi^{j+1}_{3,3}\rangle_s\\
\langle \phi^j_{2,1},\phi^{j+1}_{1,1}\rangle_s&&\langle \phi^j_{2,1},\phi^{j+1}_{2,1}\rangle_s&&....&&\langle \phi^j_{2,1},\phi^{j+1}_{3,3}\rangle_s\\
\langle \phi^j_{1,2},\phi^{j+1}_{1,1}\rangle_s&&\langle \phi^j_{1,2},\phi^{j+1}_{2,1}\rangle_s&&....&&\langle \phi^j_{1,2},\phi^{j+1}_{3,3}\rangle_s\\
\langle \phi^j_{2,2},\phi^{j+1}_{1,1}\rangle_s&&\langle \phi^j_{2,2},\phi^{j+1}_{2,1}\rangle_s&&....&&\langle \phi^j_{2,2},\phi^{j+1}_{3,3}\rangle_s\\
\end{array}\right)\left(\begin{array}{c}
b_{1,1} \\
b_{2,1}\\
b_{3,1}\\
b_{1,2}\\
b_{2,2}\\
b_{3,2}\\
b_{1,3}\\
b_{2,3}\\
b_{3,3}\\
\end{array}\right) =\left(\begin{array}{c}
0 \\
0 \\
0 \\
0 \\
0 \\
0 \\
0 \\
0 \\
0 \\
\end{array}\right).$$
Using Lemma \ref{inner}, we obtain the following equations.
$$\left(\begin{array}{ccccccccc}
 1& 1/2&-1&1/2 &2    &1/2 &-1 &1/2  &1 \\
 0&-1/2&1 &0   &-1/2 &1/2 &0  &0    &-1 \\
 0&   0&0 &0   &0    &-1/2&0  &-1/2 &1 \\
 0&   0&0 &-1/2&-1/2 & 0  &1  &1/2  &-1
\end{array}\right)\left(\begin{array}{c}
b_{1,1} \\
b_{2,1}\\
b_{3,1}\\
b_{1,2}\\
b_{2,2}\\
b_{3,2}\\
b_{1,3}\\
b_{2,3}\\
b_{3,3}
\end{array}\right)=\left(\begin{array}{c}
0 \\
0 \\
0 \\
0 \\
0 \\
0 \\
0 \\
0 \\
0 
\end{array}\right).$$

The rank of the left matrix is four, because $\phi^j_{1,1}$, $\phi^j_{2,1}$, 
$\phi^j_{1,2}$, $\phi^j_{2,2}$, are linear independent. 
So there are five solutions as shown below. 
$$\left(\begin{array}{c}
b_{1,1} \\
b_{2,1}\\
b_{3,1}\\
b_{1,2}\\
b_{2,2}\\
b_{3,2}\\
b_{1,3}\\
b_{2,3}\\
b_{3,3}
\end{array}\right)=\left(\begin{array}{c}
0 \\
0 \\
0 \\
2 \\
0 \\
0 \\
1 \\
0 \\
0 
\end{array}\right) \hbox{ or } 
\left(\begin{array}{c}
0 \\
2 \\
1 \\
0 \\
0 \\
0 \\
0 \\
0 \\
0 
\end{array}\right)  \hbox{ or  }
\left(\begin{array}{c}
1 \\
1 \\
0 \\
1 \\
-1 \\
0 \\
0 \\
0 \\
0 
\end{array}\right) \hbox{ or } 
\left(\begin{array}{c}
0 \\
0 \\
0 \\
0 \\
-1 \\
1 \\
0 \\
1 \\
1 
\end{array}\right) \hbox{ or } 
\left(\begin{array}{c}
0 \\
-1 \\
0 \\
1 \\
0 \\
-1 \\
0 \\
1 \\
0 
\end{array}\right).
$$
More precisely,  
\begin{eqnarray}
&\psi^{j,1}_{0,1}=2\phi^{j+1}_{1,2}+\phi^{j+1}_{1,3}    & \hbox{as shown in Figure 3};\\
&\psi^{j,2}_{1,0}=2\phi^{j+1}_{2,1}+\phi^{j+1}_{3,1}    & \hbox{as shown in Figure 4};\\
&\psi^{j,3}_{1,1}=-\phi^{j+1}_{2,2}+\phi^{j+1}_{3,2}+\phi^{j+1}_{2,3}+\phi^{j+1}_{3,3}
&\hbox{as shown in Figure  5};\\
&\psi^{j,4}_{1,1}=\phi^{j+1}_{1,1}+\phi^{j+1}_{2,1}+\phi^{j+1}_{1,2}-\phi^{j+1}_{2,2} 
&\hbox{as shown in Figure  6};\\
&\psi^{j,5}_{1,1}=\phi^{j+1}_{1,2}+\phi^{j+1}_{2,3}-\phi^{j+1}_{2,1}-\phi^{j+1}_{3,2} 
&\hbox{as shown in Figure  7}.
\label{FLW}
\end{eqnarray}


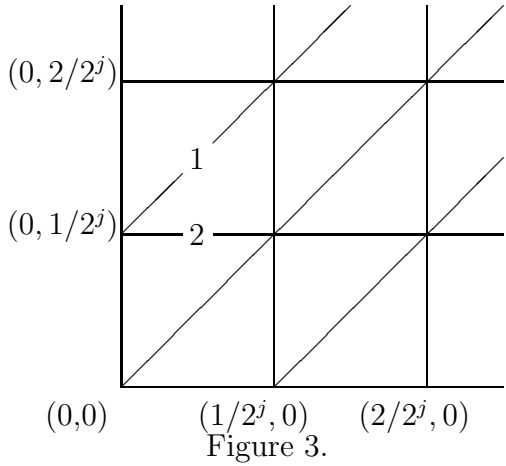
\begin{center}
\begin{picture}(7,6)(0,0)
\linethickness{0.5pt}

\put(2,1){\line(1,0){5}}
\put(2,3){\line(1,0){0.8}}
\put(3.2,3){\line(1,0){3.8}}
\put(2,5){\line(1,0){5}}

\put(2,1){\line(0,1){5}}
\put(4,1){\line(0,1){5}}
\put(6,1){\line(0,1){5}}


\put(2,3){\line(1,1){0.8}}
\put(3.2,4.2){\line(1,1){1.8}}
\put(2,1){\line(1,1){5}}
\put(4,1){\line(1,1){3}}

\put(2.89,3.86){1}
\put(2.89,2.86){2}

\put(4,0.2){\makebox(0,0){Figure 3. }}

\put(1.0,0.5){(0,0)}
\put(3.0,0.5){$(1/2^j,0)$}
\put(5.1,0.5){$(2/2^j,0)$}
\put(0.5,3){$(0,1/2^j)$}
\put(0.5,5){$(0,2/2^j)$}
\end{picture}
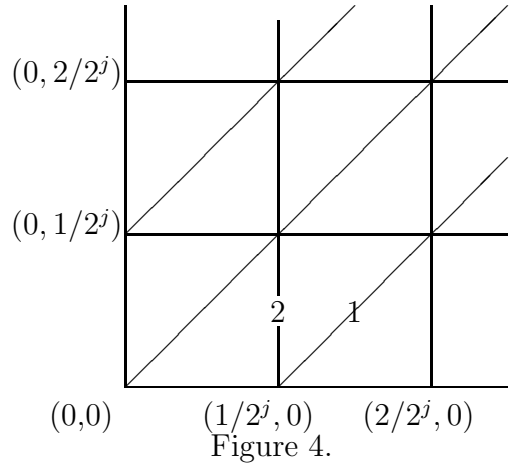
\begin{picture}(7,6)(0,0)
\linethickness{0.5pt}
\put(2,1){\line(1,0){5}}
\put(2,3){\line(1,0){5}}
\put(2,5){\line(1,0){5}}

\put(2,1){\line(0,1){5}}
\put(4,1){\line(0,1){0.8}}
\put(4,2.2){\line(0,1){3.6}}
\put(6,1){\line(0,1){5}}


\put(2,3){\line(1,1){3}}
\put(2,1){\line(1,1){5}}
\put(4,1){\line(1,1){3}}

\put(3.89,1.86){2}
\put(4.89,1.86){1}

\put(4,0.2){\makebox(0,0){Figure 4. }}

\put(1.0,0.5){(0,0)}
\put(3.0,0.5){$(1/2^j,0)$}
\put(5.1,0.5){$(2/2^j,0)$}
\put(0.5,3){$(0,1/2^j)$}
\put(0.5,5){$(0,2/2^j)$}
\end{picture}
\end{center}

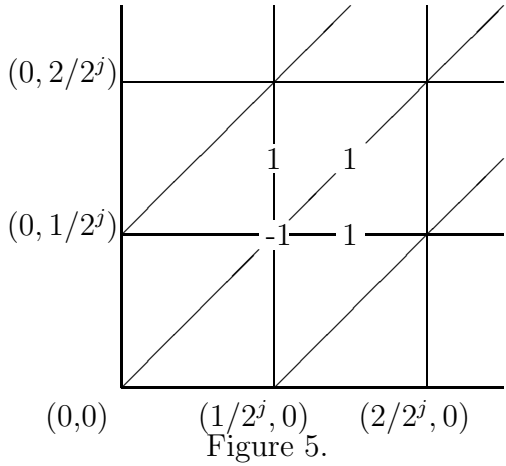
\begin{center}
\begin{picture}(7,6)(0,0)
\linethickness{0.5pt}

\put(2,1){\line(1,0){5}}
\put(2,3){\line(1,0){1.8}}
\put(4.2,3){\line(1,0){0.6}}
\put(5.2,3){\line(1,0){1.8}}
\put(2,5){\line(1,0){5}}

\put(2,1){\line(0,1){5}}
\put(4,1){\line(0,1){1.8}}
\put(4,3.2){\line(0,1){0.6}}
\put(4,4.2){\line(0,1){1.8}}
\put(6,1){\line(0,1){5}}


\put(2,3){\line(1,1){3}}
\put(2,1){\line(1,1){1.8}}
\put(4.2,3.2){\line(1,1){0.6}}
\put(5.2,4.2){\line(1,1){1.8}}
\put(4,1){\line(1,1){3}}

\put(3.89,2.86){-1}
\put(4.89,2.86){1}
\put(3.89,3.86){1}
\put(4.89,3.86){1}

\put(4,0.2){\makebox(0,0){Figure 5. }}

\put(1.0,0.5){(0,0)}
\put(3.0,0.5){$(1/2^j,0)$}
\put(5.1,0.5){$(2/2^j,0)$}
\put(0.5,3){$(0,1/2^j)$}
\put(0.5,5){$(0,2/2^j)$}
\end{picture}
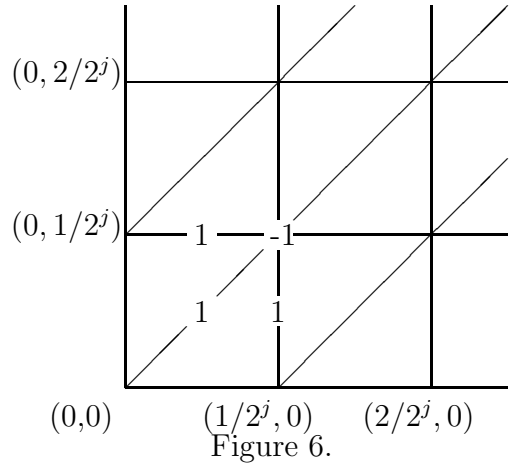
\begin{picture}(7,6)(0,0)
\linethickness{0.5pt}

\put(2,1){\line(1,0){5}}
\put(2,3){\line(1,0){0.8}}
\put(3.2,3){\line(1,0){0.6}}
\put(4.2,3){\line(1,0){2.8}}
\put(2,5){\line(1,0){5}}

\put(2,1){\line(0,1){5}}
\put(4,1){\line(0,1){0.8}}
\put(4,2.2){\line(0,1){0.6}}
\put(4,3.2){\line(0,1){2.8}}
\put(6,1){\line(0,1){5}}


\put(2,3){\line(1,1){3}}
\put(2,1){\line(1,1){0.8}}
\put(3.2,2.2){\line(1,1){0.6}}
\put(4.2,3.2){\line(1,1){2.8}}
\put(4,1){\line(1,1){3}}

\put(2.89,1.86){1}
\put(3.89,1.86){1}
\put(2.89,2.86){1}
\put(3.89,2.86){-1}

\put(4,0.2){\makebox(0,0){Figure 6. }}

\put(1.0,0.5){(0,0)}
\put(3.0,0.5){$(1/2^j,0)$}
\put(5.1,0.5){$(2/2^j,0)$}
\put(0.5,3){$(0,1/2^j)$}
\put(0.5,5){$(0,2/2^j)$}
\end{picture}
\end{center}

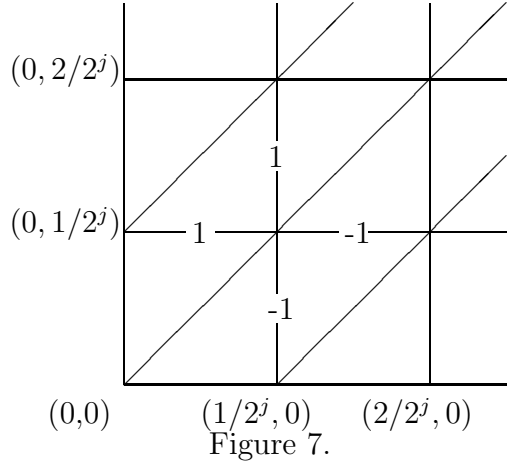
\begin{center}
\begin{picture}(7,6)(0,0)
\linethickness{0.5pt}

\put(2,1){\line(1,0){5}}
\put(2,3){\line(1,0){0.8}}
\put(3.2,3){\line(1,0){1.6}}
\put(5.2,3){\line(1,0){1.8}}
\put(2,5){\line(1,0){5}}

\put(2,1){\line(0,1){5}}
\put(4,1){\line(0,1){0.8}}
\put(4,2.2){\line(0,1){1.6}}
\put(4,4.2){\line(0,1){1.8}}
\put(6,1){\line(0,1){5}}


\put(2,3){\line(1,1){3}}
\put(2,1){\line(1,1){5}}
\put(4,1){\line(1,1){3}}

\put(2.89,2.86){1}
\put(3.89,3.86){1}
\put(3.89,1.86){-1}
\put(4.89,2.86){-1}

\put(4,0.2){\makebox(0,0){Figure 7. }}

\put(1.0,0.5){(0,0)}
\put(3.0,0.5){$(1/2^j,0)$}
\put(5.1,0.5){$(2/2^j,0)$}
\put(0.5,3){$(0,1/2^j)$}
\put(0.5,5){$(0,2/2^j)$}
\end{picture}
\end{center}

Now we consider $V_j^3$. Similarly, there are 25 non-zero coefficient for linear system 
(\ref{WWLT}) and the coefficient matrix of rank 9. 
So the dimension of solution space of $W_j^3$ is $25-9=16$.
The first five of them are the same to the wavelet functions in (7)--(\ref{FLW}). 
The other 11 are given below.

$$\begin{array}{ccccc}
&\psi^{j,1}_{0,2}=&&2\phi^{j+1}_{1,4}+\phi^{j+1}_{1,5} &\quad \quad\quad \hbox{as shown in 
Figure 8;}\\
&\psi^{j,2}_{2,0}=&&2\phi^{j+1}_{4,1}+\phi^{j+1}_{5,1} &\quad \quad\quad \hbox{as shown in
Figure 9;}\\
\quad&\psi^{j,3}_{1,2}=&&\phi^{j+1}_{3,5}+\phi^{j+1}_{3,4}+\phi^{j+1}_{2,5}-\phi^{j+1}_{2,4} 
&\quad \quad\quad \hbox{as shown in Figure 10;}\\
\quad&\psi^{j,3}_{2,2}=&&\phi^{j+1}_{5,5}+\phi^{j+1}_{5,4}+\phi^{j+1}_{4,5}
-\phi^{j+1}_{4,4} &\quad \quad\quad \hbox{as shown in Figure 11;}\\
&\psi^{j,3}_{2,1}=&&\phi^{j+1}_{5,3}+\phi^{j+1}_{5,2}+\phi^{j+1}_{4,3}-\phi^{j+1}_{4,2} 
&\quad \quad\quad \hbox{as shown in Figure 12;}\\
&\psi^{j,4}_{2,1}=&&\phi^{j+1}_{3,2}+\phi^{j+1}_{4,1}+\phi^{j+1}_{3,1}
-\phi^{j+1}_{4,2} &\quad \quad\quad \hbox{as shown in Figure 13;}\\
&\psi^{j,4}_{2,2}=&&\phi^{j+1}_{3,3}+\phi^{j+1}_{4,3}+\phi^{j+1}_{3,4}
-\phi^{j+1}_{4,4} &\quad \quad\quad \hbox{as shown in Figure 14;}\\
&\psi^{j,4}_{1,2}=&&\phi^{j+1}_{1,3}+\phi^{j+1}_{2,3}+\phi^{j+1}_{1,4}
-\phi^{j+1}_{2,4} &\quad \quad\quad \hbox{as shown in Figure 15;}\\
&\psi^{j,5}_{1,2}=&&\phi^{j+1}_{1,4}+\phi^{j+1}_{2,5}-\phi^{j+1}_{2,3}
-\phi^{j+1}_{3,4} &\quad \quad\quad \hbox{as shown in Figure 16;}\\
&\psi^{j,5}_{2,2}=&&\phi^{j+1}_{3,4}+\phi^{j+1}_{4,5}-\phi^{j+1}_{4,3}-\phi^{j+1}_{5,4} 
&\quad \quad\quad  \hbox{as shown in Figure 17;}\\
&\psi^{j,5}_{2,1}=&&\phi^{j+1}_{3,2}+\phi^{j+1}_{4,3}-\phi^{j+1}_{4,1}
-\phi^{j+1}_{5,2} &\quad \quad\quad  \hbox{as shown in Figure 18}.
\end{array}
$$

\begin{center}
\begin{picture}(7,6)(0,0)
\linethickness{0.5pt}
\put(2,2){\line(1,0){5}}
\put(2,3){\line(1,0){5}}
\put(2,4){\line(1,0){0.3}}
\put(2.7,4){\line(1,0){4.3}}
\put(2,5){\line(1,0){5}}

\put(2,2){\line(0,1){4}}
\put(3,2){\line(0,1){4}}
\put(4,2){\line(0,1){4}}
\put(5,2){\line(0,1){4}}
\put(6,2){\line(0,1){4}}


\put(2,5){\line(1,1){1}}
\put(2,4){\line(1,1){0.35}}
\put(2.7,4.7){\line(1,1){1.3}}
\put(2,3){\line(1,1){3}}
\put(2,2){\line(1,1){4}}
\put(3,2){\line(1,1){4}}
\put(4,2){\line(1,1){3}}
\put(5,2){\line(1,1){2}}

\put(2.39,4.36){1}
\put(2.39,3.86){2}

\put(4,0.5){\makebox(0,0){Figure 8. }}

\put(1.0,1.5){(0,0)}
\put(2.3,1.5){$(1/2^j,0)$}
\put(4.3,1.5){$(3/2^j,0)$}
\put(0.5,3){$(0,1/2^j)$}
\put(0.5,4){$(0,2/2^j)$}
\put(0.5,5){$(0,3/2^j)$}
\end{picture}
\begin{picture}(7,6)(0,0)
\linethickness{0.5pt}

\put(2,2){\line(1,0){5}}
\put(2,3){\line(1,0){5}}
\put(2,4){\line(1,0){5}}
\put(2,5){\line(1,0){5}}

\put(2,2){\line(0,1){4}}
\put(3,2){\line(0,1){4}}
\put(4,2){\line(0,1){0.3}}
\put(4,2.7){\line(0,1){3.3}}
\put(5,2){\line(0,1){4}}
\put(6,2){\line(0,1){4}}


\put(2,5){\line(1,1){1}}
\put(2,4){\line(1,1){2}}
\put(2,3){\line(1,1){3}}
\put(2,2){\line(1,1){4}}
\put(3,2){\line(1,1){4}}
\put(4,2){\line(1,1){0.35}}
\put(4.7,2.7){\line(1,1){2.3}}
\put(5,2){\line(1,1){2}}

\put(3.93,2.36){2}
\put(4.43,2.36){1}

\put(5,0.5){\makebox(0,0){Figure 9. }}

\put(1.0,1.5){(0,0)}
\put(2.3,1.5){$(1/2^j,0)$}
\put(4.3,1.5){$(3/2^j,0)$}
\put(0.5,3){$(0,1/2^j)$}
\put(0.5,4){$(0,2/2^j)$}
\put(0.5,5){$(0,3/2^j)$}
\end{picture}
\end{center}

\begin{center}
\begin{picture}(7,6)(0,0)
\linethickness{0.5pt}
\put(2,2){\line(1,0){5}}
\put(2,3){\line(1,0){5}}
\put(2,4){\line(1,0){0.7}}
\put(3.7,4){\line(1,0){3.3}}
\put(2,5){\line(1,0){5}}

\put(2,2){\line(0,1){4}}
\put(3,2){\line(0,1){1.7}}
\put(3,4.7){\line(0,1){1.3}}
\put(4,2){\line(0,1){4}}
\put(5,2){\line(0,1){4}}
\put(6,2){\line(0,1){4}}


\put(2,5){\line(1,1){1}}
\put(2,4){\line(1,1){2}}
\put(2,3){\line(1,1){0.7}}
\put(3.7,4.7){\line(1,1){1.3}}
\put(2,2){\line(1,1){4}}
\put(3,2){\line(1,1){4}}
\put(4,2){\line(1,1){3}}
\put(5,2){\line(1,1){2}}

\put(2.93,4.36){1}
\put(3.43,4.36){1}

\put(2.83,3.86){-1}
\put(3.43,3.86){1}

\put(5,0.5){\makebox(0,0){Figure 10. }}

\put(1.0,1.5){(0,0)}
\put(2.3,1.5){$(1/2^j,0)$}
\put(4.3,1.5){$(3/2^j,0)$}
\put(0.5,3){$(0,1/2^j)$}
\put(0.5,4){$(0,2/2^j)$}
\put(0.5,5){$(0,3/2^j)$}
\end{picture}
\begin{picture}(7,6)(0,0)
\linethickness{0.5pt}

\put(2,2){\line(1,0){5}}
\put(2,3){\line(1,0){5}}
\put(2,4){\line(1,0){1.7}}
\put(4.7,4){\line(1,0){2.3}}
\put(2,5){\line(1,0){5}}

\put(2,2){\line(0,1){4}}
\put(3,2){\line(0,1){4}}
\put(4,2){\line(0,1){1.7}}
\put(4,4.7){\line(0,1){1.3}}
\put(5,2){\line(0,1){4}}
\put(6,2){\line(0,1){4}}


\put(2,5){\line(1,1){1}}
\put(2,4){\line(1,1){2}}
\put(2,3){\line(1,1){3}}
\put(2,2){\line(1,1){1.7}}
\put(4.7,4.7){\line(1,1){1.3}}
\put(3,2){\line(1,1){4}}
\put(4,2){\line(1,1){3}}
\put(5,2){\line(1,1){2}}

\put(3.93,4.36){1}
\put(4.43,4.36){1}

\put(3.83,3.86){-1}
\put(4.43,3.86){1}

\put(5,0.5){\makebox(0,0){Figure 11. }}

\put(1.0,1.5){(0,0)}
\put(2.3,1.5){$(1/2^j,0)$}
\put(4.3,1.5){$(3/2^j,0)$}
\put(0.5,3){$(0,1/2^j)$}
\put(0.5,4){$(0,2/2^j)$}
\put(0.5,5){$(0,3/2^j)$}
\end{picture}
\end{center}

\begin{center}
\begin{picture}(7,6)(0,0)
\linethickness{0.5pt}
\put(2,2){\line(1,0){5}}
\put(2,3){\line(1,0){1.7}}
\put(4.7,3){\line(1,0){2.3}}
\put(2,4){\line(1,0){5}}
\put(2,5){\line(1,0){5}}
\put(2,2){\line(0,1){4}}
\put(3,2){\line(0,1){4}}
\put(4,2){\line(0,1){0.7}}
\put(4,3.7){\line(0,1){2.3}}
\put(5,2){\line(0,1){4}}
\put(6,2){\line(0,1){4}}


\put(2,5){\line(1,1){1}}
\put(2,4){\line(1,1){2}}
\put(2,3){\line(1,1){3}}
\put(2,2){\line(1,1){4}}
\put(3,2){\line(1,1){0.7}}
\put(4.7,3.7){\line(1,1){2.3}}
\put(4,2){\line(1,1){3}}
\put(5,2){\line(1,1){2}}

\put(3.93,3.36){1}
\put(4.43,3.36){1}

\put(3.83,2.86){-1}
\put(4.43,2.86){1}

\put(5,0.5){\makebox(0,0){Figure 12. }}

\put(1.0,1.5){(0,0)}
\put(2.3,1.5){$(1/2^j,0)$}
\put(4.3,1.5){$(3/2^j,0)$}
\put(0.5,3){$(0,1/2^j)$}
\put(0.5,4){$(0,2/2^j)$}
\put(0.5,5){$(0,3/2^j)$}
\end{picture}
\begin{picture}(7,6)(0,0)
\linethickness{0.5pt}
\put(2,2){\line(1,0){5}}
\put(2,3){\line(1,0){1.3}}
\put(4.3,3){\line(1,0){2.7}}
\put(2,4){\line(1,0){5}}
\put(2,5){\line(1,0){5}}
\put(2,2){\line(0,1){4}}
\put(3,2){\line(0,1){4}}
\put(4,2){\line(0,1){0.3}}
\put(4,3.3){\line(0,1){2.7}}
\put(5,2){\line(0,1){4}}
\put(6,2){\line(0,1){4}}
\put(2,5){\line(1,1){1}}
\put(2,4){\line(1,1){2}}
\put(2,3){\line(1,1){3}}
\put(2,2){\line(1,1){4}}
\put(3,2){\line(1,1){0.35}}
\put(4.3,3.3){\line(1,1){2.7}}
\put(4,2){\line(1,1){3}}
\put(5,2){\line(1,1){2}}
\put(3.93,2.36){1}
\put(3.43,2.36){1}
\put(3.83,2.86){-1}
\put(3.43,2.86){1}
\put(5,0.5){\makebox(0,0){Figure 13. }}
\put(1.0,1.5){(0,0)}
\put(2.3,1.5){$(1/2^j,0)$}
\put(4.3,1.5){$(3/2^j,0)$}
\put(0.5,3){$(0,1/2^j)$}
\put(0.5,4){$(0,2/2^j)$}
\put(0.5,5){$(0,3/2^j)$}
\end{picture}
\end{center}

\begin{center}
\begin{picture}(7,6)(0,0)
\linethickness{0.5pt}

\put(2,2){\line(1,0){5}}
\put(2,3){\line(1,0){5}}
\put(2,4){\line(1,0){1.3}}
\put(4.3,4){\line(1,0){2.7}}
\put(2,5){\line(1,0){5}}

\put(2,2){\line(0,1){4}}
\put(3,2){\line(0,1){4}}
\put(4,2){\line(0,1){1.3}}
\put(4,4.3){\line(0,1){1.7}}
\put(5,2){\line(0,1){4}}
\put(6,2){\line(0,1){4}}


\put(2,5){\line(1,1){1}}
\put(2,4){\line(1,1){2}}
\put(2,3){\line(1,1){3}}
\put(2,2){\line(1,1){1.35}}
\put(4.3,4.3){\line(1,1){1.7}}
\put(3,2){\line(1,1){4}}
\put(4,2){\line(1,1){3}}
\put(5,2){\line(1,1){2}}

\put(3.93,3.36){1}
\put(3.43,3.36){1}

\put(3.83,3.86){-1}
\put(3.43,3.86){1}

\put(5,0.5){\makebox(0,0){Figure 14. }}

\put(1.0,1.5){(0,0)}
\put(2.3,1.5){$(1/2^j,0)$}
\put(4.3,1.5){$(3/2^j,0)$}
\put(0.5,3){$(0,1/2^j)$}
\put(0.5,4){$(0,2/2^j)$}
\put(0.5,5){$(0,3/2^j)$}
\end{picture}
\begin{picture}(7,6)(0,0)
\linethickness{0.5pt}
\put(2,2){\line(1,0){5}}
\put(2,3){\line(1,0){5}}
\put(2,4){\line(1,0){0.3}}
\put(3.3,4){\line(1,0){3.7}}
\put(2,5){\line(1,0){5}}
\put(2,2){\line(0,1){4}}
\put(3,2){\line(0,1){1.3}}
\put(3,4.3){\line(0,1){1.7}}
\put(4,2){\line(0,1){4}}

\put(5,2){\line(0,1){4}}
\put(6,2){\line(0,1){4}}


\put(2,5){\line(1,1){1}}
\put(2,4){\line(1,1){2}}
\put(2,3){\line(1,1){0.35}}
\put(3.3,4.3){\line(1,1){1.7}}
\put(2,2){\line(1,1){4}}
\put(3,2){\line(1,1){4}}
\put(4,2){\line(1,1){3}}
\put(5,2){\line(1,1){2}}

\put(2.93,3.36){1}
\put(2.43,3.36){1}

\put(2.83,3.86){-1}
\put(2.43,3.86){1}

\put(5,0.5){\makebox(0,0){Figure 15. }}

\put(1.0,1.5){(0,0)}
\put(2.3,1.5){$(1/2^j,0)$}
\put(4.3,1.5){$(3/2^j,0)$}
\put(0.5,3){$(0,1/2^j)$}
\put(0.5,4){$(0,2/2^j)$}
\put(0.5,5){$(0,3/2^j)$}
\end{picture}
\end{center}

\begin{center}
\begin{picture}(7,6)(0,0)
\linethickness{0.5pt}
\put(2,2){\line(1,0){5}}
\put(2,3){\line(1,0){5}}
\put(2,4){\line(1,0){0.3}}
\put(3.7,4){\line(1,0){3.3}}
\put(2,5){\line(1,0){5}}
\put(2,2){\line(0,1){4}}
\put(3,2){\line(0,1){1.3}}
\put(3,4.7){\line(0,1){1.3}}
\put(4,2){\line(0,1){4}}
\put(5,2){\line(0,1){4}}
\put(6,2){\line(0,1){4}}
\put(2,5){\line(1,1){1}}
\put(2,4){\line(1,1){2}}
\put(2,3){\line(1,1){3}}
\put(2,2){\line(1,1){4}}
\put(3,2){\line(1,1){4}}
\put(4,2){\line(1,1){3}}
\put(5,2){\line(1,1){2}}
\put(2.93,4.36){1}
\put(2.83,3.36){-1}
\put(3.33,3.86){-1}
\put(2.43,3.86){1}
\put(5,0.5){\makebox(0,0){Figure 16. }}
\put(1.0,1.5){(0,0)}
\put(2.3,1.5){$(1/2^j,0)$}
\put(4.3,1.5){$(3/2^j,0)$}
\put(0.5,3){$(0,1/2^j)$}
\put(0.5,4){$(0,2/2^j)$}
\put(0.5,5){$(0,3/2^j)$}
\end{picture}
\begin{picture}(7,6)(0,0)
\linethickness{0.5pt}
\put(2,2){\line(1,0){5}}
\put(2,3){\line(1,0){5}}
\put(2,4){\line(1,0){1.3}}
\put(4.7,4){\line(1,0){2.3}}
\put(2,5){\line(1,0){5}}
\put(2,2){\line(0,1){4}}
\put(3,2){\line(0,1){4}}
\put(4,2){\line(0,1){1.3}}
\put(4,4.7){\line(0,1){1.3}}
\put(5,2){\line(0,1){4}}
\put(6,2){\line(0,1){4}}
\put(2,5){\line(1,1){1}}
\put(2,4){\line(1,1){2}}
\put(2,3){\line(1,1){3}}
\put(2,2){\line(1,1){4}}
\put(3,2){\line(1,1){4}}
\put(4,2){\line(1,1){3}}
\put(5,2){\line(1,1){2}}
\put(3.93,4.36){1}
\put(3.83,3.36){-1}
\put(4.33,3.86){-1}
\put(3.43,3.86){1}
\put(5,0.5){\makebox(0,0){Figure 17. }}
\put(1.0,1.5){(0,0)}
\put(2.3,1.5){$(1/2^j,0)$}
\put(4.3,1.5){$(3/2^j,0)$}
\put(0.5,3){$(0,1/2^j)$}
\put(0.5,4){$(0,2/2^j)$}
\put(0.5,5){$(0,3/2^j)$}
\end{picture}
\end{center}

\begin{center}
\begin{picture}(7,6)(0,0)
\linethickness{0.5pt}
\put(2,2){\line(1,0){5}}
\put(2,3){\line(1,0){1.3}}
\put(4.7,3){\line(1,0){2.3}}
\put(2,4){\line(1,0){5}}
\put(2,5){\line(1,0){5}}
\put(2,2){\line(0,1){4}}
\put(3,2){\line(0,1){4}}
\put(4,2){\line(0,1){0.3}}
\put(4,3.7){\line(0,1){2.3}}
\put(5,2){\line(0,1){4}}
\put(6,2){\line(0,1){4}}
\put(2,5){\line(1,1){1}}
\put(2,4){\line(1,1){2}}
\put(2,3){\line(1,1){3}}
\put(2,2){\line(1,1){4}}
\put(3,2){\line(1,1){4}}
\put(4,2){\line(1,1){3}}
\put(5,2){\line(1,1){2}}
\put(3.93,3.36){1}
\put(3.83,2.36){-1}
\put(4.33,2.86){-1}
\put(3.43,2.86){1}
\put(5,0.5){\makebox(0,0){Figure 18. }}
\put(1.0,1.5){(0,0)}
\put(2.3,1.5){$(1/2^j,0)$}
\put(4.3,1.5){$(3/2^j,0)$}
\put(0.5,3){$(0,1/2^j)$}
\put(0.5,4){$(0,2/2^j)$}
\put(0.5,5){$(0,3/2^j)$}
\end{picture}
\end{center}

The above computation can be carried out on $V_j^n$ for $n=3,....2^j-1$. 
We have thus obtained five types of wavelet functions: 
$$\psi^{j,1}_{0,k}=2\phi^{j+1}_{1,k+1}+\phi^{j+1}_{1,k+2} $$ 
is supported next to the vertical boundary and  is  called vertical boundary wavelet. 
$$\psi^{j,2}_{k,0}=2\phi^{j+1}_{k+1,1}+\phi^{j+1}_{k+2,1} $$ 
called horizontal boundary  wavelet, is supported next to the horizontal boundary. The
next three types are supported inside the domain. The following 
$$\psi^{j,3}_{i,k}=-\phi^{j+1}_{i+1,k+1}+\phi^{j+1}_{i+2,k+1}
+\phi^{j+1}_{i+1,k+2}+\phi^{j+1}_{i+2,k+2}$$ 
is called interior wavelet of first kind. We call 
$$\psi^{j,4}_{i,k}=-\phi^{j+1}_{2i,2k}+\phi^{j+1}_{2i-1,2k}
 +\phi^{j+1}_{2i,2k-1}+\phi^{j+1}_{2i-1,2k-1}$$ 
interior wavelet of second kind. The last one 
$$\psi^{j,5}_{i,k}=\phi^{j+1}_{2i-1,2k}+\phi^{j+1}_{2i,2k+1}-\phi^{j+1}_{2i,2k-1}
-\phi^{j+1}_{2i+1,2k}$$ 
is called interior wavelet of third kind.

\begin{theorem}
All the five types of wavelets in the $V_{j+1}^{n}$ 
are linear independent for $1\leq n \leq {2j-1}$. That is, for each $1\leq n \leq {2j-1}$,
the following functions  
$$\begin{array}{cc}
\psi^{j,1}_{0,k},    & k=1,..,n-1,\\
\psi^{j,2}_{k,0},    & k=1,..,n-1,\\
\psi^{j,3}_{i,k},& 1\leq i,k \leq n-1,\\
\psi^{j,4}_{i,k},& 1\leq i,k \leq n-1,\\
\psi^{j,5}_{i,k},& 1\leq i,k \leq n-1
\end{array}
$$
are linear independent.  
\end{theorem}
\begin{proof}
Let us prove it by induction. It is true for $n=2$ and for $n=3$. 
Suppose it is true for $n=p$, that is, 
$$\begin{array}{cc}
\psi^{j,1}_{0,k},    & k=1,..,p-1;\\
\psi^{j,2}_{k,0},    & k=1,..,p-1;\\
\psi^{j,3}_{i,k},& 1\leq i,k \leq p-1;\\
\psi^{j,4}_{i,k},& 1\leq i,k \leq p-1;\\
\psi^{j,5}_{i,k},& 1\leq i,k \leq p-1;
\end{array}
$$ are linear independent. For $n=p+1$, there are $6p-1$ new functions which are 
$$
\begin{array}{cc}
\psi^{j,1}_{0,k},    &k=p;\\
\psi^{j,2}_{k,0},    &k=p;\\
\psi^{j,3}_{i,k},& i \quad or\quad k =p;\\
\psi^{j,4}_{i,k},& i \quad or\quad k =p;\\
\psi^{j,5}_{i,k},& i\quad or\quad k =p. 
\end{array}
$$
Suppose they are not linear independent.  That is, one can find 
$$
\begin{array}{cc}
a^1,  &  \\
a^2,  &  \\
a^3_{i,k},& i \quad or\quad k =p;\\
a^4_{i,k},& i \quad or\quad k =p;\\
a^5_{i,k},& i\quad or\quad k =p\\
\end{array}
$$
such that 
\begin{equation}
a^1\psi^{j,1}_{0,p}+a^2\psi^{j,2}_{p,0}+\sum_{i \hbox{ or } k =p}
a^3_{i,k}\psi^{j,3}_{i,k}
+\sum_{i \hbox{ or } k =p}a^4_{i,k}\psi^{j,4}_{i,k}
+\sum_{i \hbox{ or } k =p}a^5_{i,k}\psi^{j,5}_{i,k}+\psi^{'}=0,
\label{equation7}
\end{equation}
where $\psi^{'}$ is linear combination of the following functions:
$$\begin{array}{cc}
\psi^{j,1}_{0,k},    & k=1,..,p-1;\\
\psi^{j,2}_{k,0},    & k=1,..,p-1;\\
\psi^{j,3}_{i,k},& 1\leq i,k \leq p-1;\\
\psi^{j,4}_{i,k},& 1\leq i,k \leq p-1;\\
\psi^{j,5}_{i,k},& 1\leq i,k \leq p-1.
\end{array}
$$
By the definition, $\phi^{j+1}_{2i+1,2k+1}, i=p$ or $k=p$  appear only once 
in $\psi^{j,3}_{i,k}, i=p$ or $k=p$ , $\psi^{j,1}_{0,p}$ and $\psi^{j,2}_{p,0}$. 
Since $\phi^{j+1}$ are linear independent, that is,  
$a^3_{i,k}=0,i \quad or\quad k =p$, $a^1=0$, and $a^2=0$. 
Thus the equation (\ref{equation7}) can be simplified to  
\begin{equation}
\sum_{i \, or\, k =p}a^4_{i,k}\psi^{j,4}_{i,k}+\sum_{i \, or\, k =p}
a^5_{i,k}\psi^{j,5}_{i,k}+\psi^{'}=0.
\label{equation8}
\end{equation}
By the similar reason, $\phi^{j+1}_{2i,2k}, i=p$ or $k=p$  appear only once 
in   $\psi^{j,4}_{i,k}, i=p$ or $k=p$. Since $\phi^{j+1}_{ik}$ are linear independent, 
\(a^4_{i,k}=0, i \, or \, k =p\). Thus the equation (\ref{equation8}) 
can be further simplified to the following equation 
\begin{equation*}
\sum_{i \, or\, k =p}a^5_{i,k}\psi^{j,5}_{i,k}+\psi^{'}=0.
\end{equation*}
Similarly, \(a^5_{i,k}=0, i \, or \, k =p\) too. Thus the equation (\ref{equation7}) 
is reduced to  
\begin{equation*}
\psi^{'}=0.
\end{equation*}
By induction hypothesis, all the coefficient of $\psi^{'}=0$ are zeros. Hence, 
$$\begin{array}{cc}
\psi^{j,1}_{0,k},    & k=1,..,n-1,\\
\psi^{j,2}_{k,0},    & k=1,..,n-1,\\
\psi^{j,3}_{i,k},& 1\leq i,k \leq n-1,\\
\psi^{j,4}_{i,k},& 1\leq i,k \leq n-1,\\
\psi^{j,5}_{i,k},& 1\leq i,k \leq n-1
\end{array}
$$
are linear independent.
\end{proof}

\begin{theorem}
All the five types of wavelets in the $W_j^{n}$ form a  
basis of $W_j^{n}$ for  $1\leq n \leq {2j-1}$. That is, 
$$W_j^{n}=\hbox{span} \{
\psi^{j,1}_{0,k},  \psi^{j,2}_{k,0},   
\psi^{j,3}_{i,k}, 
\psi^{j,4}_{i,k}, 
\psi^{j,5}_{i,k},  1\leq i,k \leq n-1\}
$$
for $1\leq n \leq {2j-1}$.
\end{theorem}
\begin{proof}
The dimension of $W_j^{n}$ is ${(2n-1)}^2-{(n)}^2=3n^2-4n+1$. It is easy to count that 
there are ${(2n-1)}^2-{(n)}^2=3n^2-4n+1$ functions in the following set 
$$
\begin{array}{cc}
\psi^{j,1}_{0,k},    & k=1,..,n;\\
\psi^{j,2}_{k,0},    & k=1,..,n;\\
\psi^{j,3}_{i,k},& 1\leq i,k \leq n;\\
\psi^{j,4}_{i,k},& 1\leq i,k \leq n;\\
\psi^{j,5}_{i,k},& 1\leq i,k \leq n
\end{array}
$$
which  all belong to the space $W_j^{n}$. Since they are linear independent, 
they form a basis for space  $W_j^{n}$,  where $1\leq n \leq {2j-1}$.
\end{proof}

Finally we need to find wavelets in $W_j^{2^j}\backslash W_j^{2^j-1}$. The computations  
are the same to the above except for that there is  one 
globally supported basis function. In fact  
the following pictures show the basis functions located on the top boundary of the domain 
$\Omega$. (We omit the pictures  for the basis functions on the right vertical boundary 
which are symmetric with respect to the line y=x are those  basic functions on the 
top horizontal boundary of $\Omega$.)

\begin{picture}(7,6)(0,0)
\linethickness{0.5pt}
\put(2,3){\line(1,0){5}}
\put(2,4){\line(1,0){5}}
\put(2,5){\line(1,0){5}}

\put(2,2){\line(0,1){3}}
\put(3,2){\line(0,1){3}}
\put(4,2){\line(0,1){3}}
\put(5,2){\line(0,1){3}}
\put(6,2){\line(0,1){3}}

\put(2.5,5.3){($\frac{1}{2^j}$,1)}
\put(4.5,5.3){($\frac{3}{2^j}$,1)}
\put(0.7,3){(0,$\frac{2^j-2}{2^j}$)}
\put(0.7,4){(0,$\frac{2^j-1}{2^j}$)}
\put(0.7,5){(0,$\frac{2^j}{2^j}$)}

\put(2,4){\line(1,1){1}}
\put(2,3){\line(1,1){2}}
\put(2,2){\line(1,1){3}}
\put(3,2){\line(1,1){3}}
\put(4,2){\line(1,1){3}}
\put(5,2){\line(1,1){2}}

\put(2.39,4.37){1}
\put(2.39,3.87){2}
\put(4,1.5){\makebox(0,0){Figure 19. }}
\end{picture}
\begin{picture}(7,5)(0,0)
\linethickness{0.5pt}

\put(2,3){\line(1,0){5}}
\put(2,4){\line(1,0){5}}
\put(2,5){\line(1,0){5}}

\put(2,2){\line(0,1){3}}
\put(3,2){\line(0,1){3}}
\put(4,2){\line(0,1){3}}
\put(5,2){\line(0,1){3}}
\put(6,2){\line(0,1){3}}

\put(2.5,5.3){($\frac{1}{2^j}$,1)}
\put(4.5,5.3){($\frac{3}{2^j}$,1)}
\put(0.7,3){(0,$\frac{2^j-2}{2^j}$)}
\put(0.7,4){(0,$\frac{2^j-1}{2^j}$)}
\put(0.7,5){(0,$\frac{2^j}{2^j}$)}

\put(2,4){\line(1,1){1}}
\put(2,3){\line(1,1){2}}
\put(2,2){\line(1,1){3}}
\put(3,2){\line(1,1){3}}
\put(4,2){\line(1,1){3}}
\put(5,2){\line(1,1){2}}

\put(3.02,4.37){1}
\put(2.89,3.87){-1}
\put(3.39,4.37){1}
\put(3.39,3.87){1}


\put(4,1.5){\makebox(0,0){Figure 20. }}
\end{picture}

\begin{picture}(7,6)(0,0)
\linethickness{0.5pt}

\put(2,3){\line(1,0){5}}
\put(2,4){\line(1,0){5}}
\put(2,5){\line(1,0){5}}

\put(2,2){\line(0,1){3}}
\put(3,2){\line(0,1){3}}
\put(4,2){\line(0,1){3}}
\put(5,2){\line(0,1){3}}
\put(6,2){\line(0,1){3}}

\put(2.5,5.3){($\frac{1}{2^j}$,1)}
\put(4.5,5.3){($\frac{3}{2^j}$,1)}
\put(0.7,3){(0,$\frac{2^j-2}{2^j}$)}
\put(0.7,4){(0,$\frac{2^j-1}{2^j}$)}
\put(0.7,5){(0,$\frac{2^j}{2^j}$)}

\put(2,4){\line(1,1){1}}
\put(2,3){\line(1,1){2}}
\put(2,2){\line(1,1){3}}
\put(3,2){\line(1,1){3}}
\put(4,2){\line(1,1){3}}
\put(5,2){\line(1,1){2}}

\put(3.89,3.84){-1}
\put(4.03,3.37){1}
\put(3.39,3.87){1}
\put(3.39,3.37){1}
\put(4,1.5){\makebox(0,0){Figure 21. }}
\end{picture}
\begin{picture}(7,5)(0,0)
\linethickness{0.5pt}

\put(2,3){\line(1,0){5}}
\put(2,4){\line(1,0){5}}
\put(2,5){\line(1,0){5}}

\put(2,2){\line(0,1){3}}
\put(3,2){\line(0,1){3}}
\put(4,2){\line(0,1){3}}
\put(5,2){\line(0,1){3}}
\put(6,2){\line(0,1){3}}

\put(2.5,5.3){($\frac{1}{2^j}$,1)}
\put(4.5,5.3){($\frac{3}{2^j}$,1)}
\put(0.7,3){(0,$\frac{2^j-2}{2^j}$)}
\put(0.7,4){(0,$\frac{2^j-1}{2^j}$)}
\put(0.7,5){(0,$\frac{2^j}{2^j}$)}

\put(2,4){\line(1,1){1}}
\put(2,3){\line(1,1){2}}
\put(2,2){\line(1,1){3}}
\put(3,2){\line(1,1){3}}
\put(4,2){\line(1,1){3}}
\put(5,2){\line(1,1){2}}

\put(2.84,4.37){1}
\put(2.89,3.37){-1}
\put(3.39,3.84){-1}
\put(2.39,3.87){1}
\put(4,1.5){\makebox(0,0){Figure 22. }}
\end{picture}

\begin{picture}(7,5)(0,0)
\linethickness{0.5pt}

\put(2,3){\line(1,0){5}}
\put(2,4){\line(1,0){5}}
\put(2,5){\line(1,0){5}}

\put(2,2){\line(0,1){3}}
\put(3,2){\line(0,1){3}}
\put(4,2){\line(0,1){3}}
\put(5,2){\line(0,1){3}}
\put(6,2){\line(0,1){3}}

\put(2.5,5.3){($\frac{1}{2^j}$,1)}
\put(4.5,5.3){($\frac{3}{2^j}$,1)}
\put(0.7,3){(0,$\frac{2^j-2}{2^j}$)}
\put(0.7,4){(0,$\frac{2^j-1}{2^j}$)}
\put(0.7,5){(0,$\frac{2^j}{2^j}$)}

\put(2,4){\line(1,1){1}}
\put(2,3){\line(1,1){2}}
\put(2,2){\line(1,1){3}}
\put(3,2){\line(1,1){3}}
\put(4,2){\line(1,1){3}}
\put(5,2){\line(1,1){2}}

\put(2.39,4.37){1}
\put(2.89,4.37){2}

\put(6,1.5){\makebox(0,0){Figure 23. }}
\end{picture}

\begin{picture}(15,5)(0,0)
\linethickness{0.5pt}

\put(2,3){\line(1,0){7}}
\put(2,4){\line(1,0){7}}
\put(2,5){\line(1,0){7}}

\put(2,2){\line(0,1){3}}
\put(3,2){\line(0,1){3}}
\put(4,2){\line(0,1){3}}
\put(5,2){\line(0,1){3}}
\put(6,2){\line(0,1){3}}
\put(7,2){\line(0,1){3}}
\put(8,2){\line(0,1){3}}
\put(9,2){\line(0,1){3}}
\put(10,2){\line(0,1){3}}
\put(11,2){\line(0,1){3}}


\put(9,3){\line(1,0){4}}
\put(9,4){\line(1,0){4}}
\put(9,5){\line(1,0){4}}

\put(10,2){\line(0,1){3}}
\put(11,2){\line(0,1){3}}
\put(12,2){\line(0,1){3}}
\put(13,2){\line(0,1){3}}

\put(10.3,5.3){($\frac{2^j-2}{2^j}$,1)}
\put(2.5,5.3){($\frac{1}{2^j}$,1)}
\put(12.6,5.3){(1,1)}
\put(0.7,3){(0,$\frac{2^j-2}{2^j}$)}
\put(0.7,4){(0,$\frac{2^j-1}{2^j}$)}
\put(0.7,5){(0,$\frac{2^j}{2^j}$)}

\put(2,4){\line(1,1){1}}
\put(2,3){\line(1,1){2}}
\put(2,2){\line(1,1){3}}
\put(3,2){\line(1,1){3}}
\put(4,2){\line(1,1){3}}
\put(5,2){\line(1,1){3}}
\put(6,2){\line(1,1){3}}
\put(7,2){\line(1,1){3}}
\put(8,2){\line(1,1){3}}
\put(9,2){\line(1,1){3}}
\put(10,2){\line(1,1){3}}

\put(9,4){\line(1,1){1}}
\put(9,3){\line(1,1){2}}
\put(9,2){\line(1,1){3}}
\put(10,2){\line(1,1){3}}
\put(11,2){\line(1,1){2}}
\put(12,2){\line(1,1){1}}


\put(2.39,4.37){1}
\put(3.39,4.37){1}
\put(4.39,4.37){1}
\put(5.39,4.37){1}
\put(6.39,4.37){1}
\put(7.39,4.37){1}
\put(8.39,4.37){1}
\put(9.39,4.37){1}
\put(10.39,4.37){1}
\put(11.39,4.37){1}
\put(12.39,4.37){1}
\put(6,1.5){\makebox(0,0){Figure 24. }}
\end{picture} 
The last one (cf. Figure 24) is the only special basis function 
since it is not local supported.  
The numbers of all these wavelets in $W_j^{2^j}\backslash W_j^{2^j-1}$ 
amount to ${2^{j+3}}-8$ which is equal to the number of dimension 
of $V_{j+1}^{2^j}\backslash V_{j+1}^{2^j-1}$.

\begin{theorem}
All the wavelets in the $W_j^{2^j}\backslash W_j^{2^j-1}$ are linear independent 
and form a basis for $V_{j+1}^{2^j}\backslash V_{j+1}^{2^j-1}$ which is spanned
by the functions in $\{\phi^{j+1}_{i,k}, 2^{j+1}-2\leq i,k \leq 2^{j+1}-1\}$.  
\end{theorem}
\begin{proof}
Let us just concentrate on the basis functions in 
$V_{j+1}^{2^j}\backslash V_{j+1}^{2^j-1}$ and 
in  $W_j^{2^j}\backslash W_j^{2^j-1}$. 
Then the scaling matrix between two sets of basis functions is 
the following matrix up to a constant

$$A=\left(\begin{array}{cccccccccccc}
D\\
B1&B2\\
&B1&B2\\
&&B1&B2\\
&&&\ddots\\
&&&&B1&B2\\
&&&&&C1&C2\\
C3&C3&C3&\dots&C3&C3&C3\\
&&&&&&C4\\
&&&&&&B2'&B1'\\
&&&&&&&B2'&B1'\\
&&&&&&&&\ddots\\
&&&&&&&&&B2'&B1'\\
&&&&&&&&&&D'\\
\end{array}\right),
$$

where 
$$
D=\left(\begin{array}{cccc}
1&2&0&0\\
\end{array}\right),
\quad
B1=\left(\begin{array}{cccc}
1&0&2\\
&1&0&-1\\
&1&1&0\\
&&1&-1
\end{array}\right),
\quad
B2=\left(\begin{array}{cccccccccc}
0&0&0&0\\
0&0&0&0\\
0&-1&0&0\\
1&1&0&0
\end{array}\right),
$$

$$
D'=\left(\begin{array}{cccc}
0&0&2&1\\
\end{array}\right),
\quad
B1'=\left(\begin{array}{cccc}
-1&1\\
0&1&1\\
-1&0&1\\
&2&0&1
\end{array}\right),
\quad
B2'=\left(\begin{array}{cccc}
0&0&1&1\\
0&0&-1&0\\
0&0&0&0\\
0&0&0&0
\end{array}\right),
$$
$$
C1=\left(\begin{array}{cccc}
1&0&2\\
&1&0&-1\\
&1&1&0\\
&&1&-1\\
\end{array}\right),
\quad
C2=\left(\begin{array}{ccccc}
0&0&0&0\\
0&0&1&0\\
0&-1&-1&0\\
1&1&0&0\\
\end{array}\right),
\quad
C4=\left(\begin{array}{ccccc}
0&2&0&1
\end{array}\right),
$$
$$C3=\left(\begin{array}{cccc}
1&0&0&0
\end{array}\right).
$$
Let $E=(m \,\, n \,\, 0 \,\, 0)$. By the row operations we have 
$$\left(\begin{array}{cccccccccc}
E\\
B1&B2\\
&B1&B2
\end{array}\right)=\left(\begin{array}{cccccccccc}
m&n&0\\
1&0&2\\
&1&0&-1\\
&1&1&0&0&-1\\
&&1&-1&1&1\\
&&&&1&0&2\\
&&&&&1&0&-1\\
&&&&&1&1&0&0&-1\\
&&&&&&1&-1&1&1\\
\end{array}\right)$$
$$\rightarrow\left(\begin{array}{cccccccccc}
m&n&&&&&&&&\\
&-n&2m&&&&&&&\\
&&2m&-n&&&&&&\\
&&&n&m&&&&&\\
&&&&2m+n&2n&0&0&\\
&&&&1&0&2\\
&&&&&1&0&-1\\
&&&&&1&1&0&0&-1\\
&&&&&&1&-1&1&1\\
\end{array}\right).
$$
Similar  for $B'$. Thus by row operations, 
$$A\rightarrow\left(\begin{array}{ccccccccccc}
A_1&G_1&&&&&&&&\\
&A_2&G_2&&&&&&&\\
&&A_3&G_3&&&&&&\\
&&&\ddots&&&&&&\\
&&&&A_{2^j-2}&G_{2^j-2}&&&&\\
&&&&&C'_1&C'_2&&&\\
&&&&&&G'_{2^j-2}&A'_{2^j-2}&&&\\
&&&&&&&&\ddots&&\\
&&&&&&&&G'_2&A'_2&\\
&&&&&&&&&G'_1&A'_1\\
\end{array}\right),$$
where $A_n$ is an upper triangular matrix of size $4\times 4$ while  
$A_n'$ is a lower triangular  matrix of size $4\times 4$ which are given below.

\[A_1=\left(\begin{array}{cccc}
1&2&0&0\\
&-1&1&0\\
&&1&-1\\
&&&2
\end{array}\right), \,
G_1=\left(\begin{array}{cccc}
0&0&0&0\\
0&0&0&0\\
0&0&0&0\\
1&0&0&0\\
\end{array}\right),\,
A_2=\left(\begin{array}{cccc}
1&1&0&0\\
&-1&2&0\\
&&2&-1\\
&&&1
\end{array}\right),\]
\[
G_2=\left(\begin{array}{cccc}
0&0&0&0\\
0&0&0&0\\
0&0&0&0\\
1&0&0&0\\
\end{array}\right), \,
A_n=\left(\begin{array}{cccc}
n&2&0&0\\
&-1&n&0\\
&&n&-1\\
&&&2
\end{array}\right), \,
G_n=\left(\begin{array}{cccc}
0&0&0&0\\
0&0&0&0\\
0&0&0&0\\
n&0&0&0\\
\end{array}\right),
\]
\[
A'_n=\left(\begin{array}{cccc}
2&&&\\
-1&n&&\\
0&n&-1&\\
0&0&2&n\\
\end{array}\right), \,
G'_n=\left(\begin{array}{cccc}
0&0&0&n\\
0&0&0&0\\
0&0&0&0\\
0&0&0&0\\
\end{array}\right)
\]

and the matrix $(C'_1 \quad C'_2)$ is the following matrix

$$(C'_1 \quad C'_2)=\left(\begin{array}{ccccccccc}
2^{j+1}-5&2\\
1&0&2&\\
&1&0&-1&0&0&1\\
&1&1&0&0&-1&-1\\
&&1&-1&1&1&\\
2^{j+1}-5&0&0&0&1&0&\\
&&&&&2&0&1\\
&&&&&&2&2^{j+1}-5
\end{array}\right).
$$
It is easy to see the rank of $(C'_1 \quad C'_2)$ is 8. 
Thus the rank of A is $8(2^j)-8$. Thus, all the prewavelet functions constructed above  
in the $W_j^{2^j}\backslash W_j^{2^j-1}$ are linear independent and hence
form a basis  of $V_{j+1}^{2^j}\backslash V_{j+1}^{2^j-1}$.
\end{proof}

It is easy to see that the coefficients of the prewavelet functions in $ W_j^{2^j-1}$ 
in terms of the basis functions of $V_{j+1}^{2^j}\backslash V_{j+1}^{2^j-1}$ are all zeros. 
Thus the prewavelet functions in $W_j^{2^j-1}$ together with the prewavelet functions
in  $V_{j+1}^{2^j}\backslash V_{j+1}^{2^j-1}$ are linear independent. 
It follows the main result in this paper.
\begin{theorem}
All the prewavelet functions in the $W_j^{2^j}\backslash W_j^{2^j-1}$ 
and the prewavelet functions in $W_j^{2^j-1}$ form a basis for $W_j$.
\end{theorem}

\section{The Prewavelet Method for Poisson Equation}
Let us use the basis functions of $V_j$ and $W_j$ to solve Poisson equation (\ref{PDE}).
Mainly we explain how to compute $h_j\in W_j$. Let $g_j\in V_j$ and $g_{j+1}\in V_{j+1}$
be two FEM solutions. We aim to show that $h_j+g_j= g_{j+1}$.
 
By a reordering the indices $(i,k), 1\le i, k\le 2^j$ in a linear fashion, 
let $V_j=span\{\phi^j_1,......,\phi^j_{N_j}\}$. Also, we reorder all five type 
wavelet functions as well as the globally supported wavelet to denote 
$W_j=span\{\psi^j_1,......\psi^j_{N_{j+1}-N_j}\}$. 
Let $\Phi^j$, $\Psi^j$ be following vectors,
$$\Phi^j=\left(\begin{array}{c}
\phi^j_{1} \\
\phi^j_{2} \\
\vdots\\
\phi^j_{N_j}
\end{array}\right),\quad\quad\quad
\Psi^j=\left(\begin{array}{c}
\psi^j_{1} \\
\psi^j_{2} \\
\vdots\\
\psi^j_{N_{j+1}-N_j}
\end{array}\right).
$$
Then we have the following equations
\begin{equation*}
\begin{split}
\Phi^j=B_{j}\Phi^{j+1}, \quad\quad\quad \Psi^j=C_j\Phi^{j+1},
\end{split}
\end{equation*}
where $B_{j}$ is $N_j\times N_{j+1}$ refinable matrix, 
and $C_j$ is a wavelet matrix of size $(N_{j+1}-N_j)\times N_{j+1}$.
Let $D_j$ and $E_j$ be the following matrices:
$$D_j=\left(\begin{array}{ccccccc}
\langle\phi^j_{1},\phi^{j}_{1}\rangle_s&&\langle\phi^j_{1},\phi^{j}_{2}\rangle_s
&&\cdots\cdots&&\langle\phi^j_{1},\phi^{j}_{N_{j}}\rangle_s \\
\langle\phi^j_{2},\phi^{j}_{1}\rangle_s&&\langle\phi^j_{2},\phi^{j}_{2}\rangle_s
&&\cdots\cdots&&\langle\phi^j_{2},\phi^{j}_{N_{j}}\rangle_s \\
\vdots&&\vdots&&\ddots&&\vdots\\
\langle\phi^j_{N_j},\phi^{j}_{1}\rangle_s&&\langle\phi^j_{N_j},\phi^{j}_{2}\rangle_s
&&\cdots\cdots&&\langle\phi^j_{N_j},\phi^{j}_{N_{j}}\rangle_s
\end{array}\right)
$$
$$E_j=\left(\begin{array}{ccccccc}
\langle\psi^j_{1},\psi^{j}_{1}\rangle_s&&\langle\psi^j_{1},\psi^{j}_{2}\rangle_s
&&\cdots\cdots&&\langle\psi^j_{1},\psi^{j}_{N_{j+1}-N_j}\rangle_s \\
\langle\psi^j_{2},\psi^{j}_{1}\rangle_s&&\langle\psi^j_{2},\psi^{j}_{2}\rangle_s
&&\cdots\cdots&&\langle\psi^j_{2},\psi^{j}_{N_{j+1}-N_j}\rangle_s \\
\vdots&&\vdots&&\ddots&&\vdots\\
\langle\psi^j_{N_{j+1}-N_j},\psi^{j}_{1}\rangle_s&&\langle\psi^j_{N_{j+1}-N_j},
\psi^{j}_{2}\rangle_s&&\cdots\cdots&&\langle\psi^j_{N_{j+1}-N_j},\psi^{j}_{N_{j+1}-N_j}
\rangle_s. 
\end{array}\right).
$$
It is easy to see that $B_{j}D_{j+1}C_j^T=0$ is equivalent to 
$V_j \bot W_j$. Clearly, we have $D_j=B_{j}D_{j+1}B_j^T$ and $E_j=C_jD_{j+1}C_j^T$.

Let $g_j$ be  the projection  of g in $V_j$, and $h_j$ be the projection of g in $W_j$. 
Since $V_j\bigoplus W_j=V_{j+1}$,  $g_j+h_j$ will be equal to $g_{j+1}$. Let us write 
$g_j=\sum_{j=1}^{N_j}a_i \phi^j_i=(a_1,a_2,....,a_{N_j})\Phi^j$. 
Similarly, $h_j=(b_1,b_2,....,b_{N_{j+1}-N_j})\Psi^j$, and 
$g_{j+1}=(c_1,c_2,....,c_{N_{j+1}})\Phi^{j+1}$. By computing the weak solutions
$h_j, g_j$, and $g_{j+1}$ in $W_j, V_j$, and $V_{j+1}$, respectively, we have
$$
D_j\left(\begin{array}{c}
a_{1}\\
a_{2}\\ 
\vdots\\
a_{N_{j}}
\end{array}\right)=\left(\begin{array}{c}
\langle\phi^{j}_{1},g\rangle\\
\langle\phi^{j}_{2},g\rangle\\
\vdots\\
\langle\phi^{j}_{N_{j}},g\rangle\\
\end{array}\right),
$$
$$
E_j\left(\begin{array}{c}
b_{1}\\
b_{2}\\ 
\vdots\\
b_{N_{j+1}-N_{j}}
\end{array}\right)=\left(\begin{array}{c}
\langle\psi^{j}_{1},g\rangle\\
\langle\psi^{j}_{2},g\rangle\\
\vdots\\
\langle\psi^{j}_{N_{j+1}-N_{j}}, g\rangle
\end{array}\right),
$$
$$
D_{j+1}\left(\begin{array}{c}
c_{1}\\
c_{2}\\ 
\vdots\\
c_{N_{j+1}}
\end{array}\right)=\left(\begin{array}{c}
\langle\phi^{j+1}_{1}, g\rangle\\
\langle\phi^{j+1}_{2}, g\rangle\\
\vdots\\
\langle\phi^{j+1}_{N_{j+1}}, g\rangle
\end{array}\right).
$$

It follows
$$
\left(\begin{array}{c}
a_{1}\\
a_{2}\\ 
\vdots\\
a_{N_{j}}
\end{array}\right)={(D_j)}^{-1}\left(\begin{array}{c}
\langle\phi^{j}_{1}, g\rangle\\
\langle\phi^{j}_{2}, g\rangle\\
\vdots\\
\langle\phi^{j}_{N_{j}}, g\rangle
\end{array}\right),
$$
$$
\left(\begin{array}{c}
b_{1}\\
b_{2}\\ 
\vdots\\
b_{N_{j+1}-N_{j}}
\end{array}\right)={(E_j)}^{-1}\left(\begin{array}{c}
\langle\psi^{j}_{1}, g\rangle\\
\langle\psi^{j}_{2}, g\rangle\\
\vdots\\
\langle\psi^{j}_{N_{j+1}-N_{j}}, g\rangle
\end{array}\right),
$$
$$
\left(\begin{array}{c}
c_{1}\\
c_{2}\\ 
\vdots\\
c_{N_{j+1}}
\end{array}\right)={(D_{j+1})}^{-1}\left(\begin{array}{c}
\langle\phi^{j+1}_{1}, g\rangle\\
\langle\phi^{j+1}_{2}, g\rangle\\
\vdots\\
\langle\phi^{j+1}_{N_{j+1}}, g\rangle
\end{array}\right).
$$
The above linear systems provide a computational method to find  $g_j$, $h_j$. 

We now show $h_j+g_j=g_{j+1}$. That is, $g_{j+1}$ can be computed 
by using $h_j$ and $g_j$ only. Indeed, we have 
\begin{equation*}
\begin{split}
g_j&=(a_1,a_2,......,a_{N_j})\Phi^j =({\Phi^j})^T{(a_1,a_2,......,a_{n_j})}^T\\
   &=({\Phi^{j+1}})^TB_{j}^T{(a_1,a_2,......,a_{N_j})}^T \\
   &=({\Phi^{j+1}})^TB_{j}^TD_j^{-1}{(\langle\phi^{j}_{1},g\rangle,
\langle\phi^{j}_{2}, g\rangle,\cdots\langle\phi^{j}_{N_{j}},g\rangle)}^T \\  
   &=({(\Phi^{j+1})})^TB_{j}^T{(B_{j}D_{j+1}B_{j}^T)}^{-1}B_{j}
{(\langle\phi^{j+1}_{1},g\rangle,\langle\phi^{j+1}_{2}, g\rangle,
\cdots\langle\phi^{j+1}_{N_{j+1}}, g\rangle)}^T .\\  
\end{split}
\end{equation*}

Similarly,
$$\begin{array}{ccc}
h_j&=&({(\Phi^{j+1})})^TC_{j}^T{(C_{j}D_{j+1}C_{j}^T)}^{-1}C_{j}
{(\langle\phi^{j+1}_{1},g\rangle,\langle\phi^{j+1}_{2},g\rangle,
\cdots\langle\phi^{j+1}_{N_{j+1}},g\rangle)}^T .\\  
\end{array}$$
and
$$\begin{array}{ccc}
g_{j+1}&=&({(\Phi^{j+1})})^TD_{j+1}^{-1}{(\langle\phi^{j+1}_{1}, g\rangle,
\langle\phi^{j+1}_{2}, g\rangle,\cdots\langle\phi^{j+1}_{N_{j+1}}, g\rangle)}^T .\\  
\end{array}$$
In order to show $h_j+g_j=g_{j+1}$,  we only need to prove 
\begin{equation}
\label{key}
B_{j}^T{(B_{j}D_{j+1}B_{j}^T)}^{-1}B_{j}+C_{j}^T{(C_{j}D_{j+1}
C_{j}^T)}^{-1}C_{j}=D_{j+1}^{-1}.
\end{equation}
Notice that $B_{j}$ and $C_{j}$ are not square matrices.  That is we can not invert 
$B_j$ and $C_j$. Consider 
\[\left(\begin{array}{c}
B_{j} \\
C_{j} 
\end{array}\right)
D_{j+1}
\left(\begin{array}{cc}
B_{j}^T&C_{j}^T
\end{array}\right)=
\left(\begin{array}{ccc}
B_{j}D_{j+1}B_{j}^T&&B_{j}D_{j+1}C_{j}^T \\
C_{j}D_{j+1}B_{j}^T&&C_{j}D_{j+1}C_{j}^T 
\end{array}\right)\]
\[=\left(\begin{array}{ccc}
B_{j}D_{j+1}B_{j}^T&&0 \\
0&&C_{j}D_{j+1}C_{j}^T 
\end{array}\right)
\]
by using the orthogonal conditions of $V_j$ and $W_j$. 
Then we have the following equation
$$\left(\begin{array}{c}
B_{j}D_{j+1} \\
C_{j}D_{j+1} 
\end{array}\right)
\left(\begin{array}{cc}
B_{j}^T& C_{j}^T
\end{array}\right)
\left(\begin{array}{ccc}
{(B_{j}D_{j+1}B_{j}^T)}^{-1}&&0 \\
0&&{(C_{j}D_{j+1}C_{j}^T)}^{-1} 
\end{array}\right)=I,
$$
where $I$ stands for the identity matrix. In other words, we have 
$$\left(\begin{array}{c}
B_{j}D_{j+1} \\
C_{j}D_{j+1} 
\end{array}\right)
\left(\begin{array}{cc}
B_{j}^T{(B_{j}D_{j+1}B_{j}^T)}^{-1}&C_{j}^T{(C_{j}D_{j+1}C_{j}^T)}^{-1}
\end{array}\right)
=I
$$
which can be rewritten in the following form
$$\left(\begin{array}{cccc}
B_{j}^T{(B_{j}D_{j+1}B_{j}^T)}^{-1}&C_{j}^T{(C_{j}D_{j+1}C_{j}^T)}^{-1}
\end{array}\right)\left(\begin{array}{c}
B_{j}D_{j+1} \\
C_{j}D_{j+1} 
\end{array}\right)
=I.
$$
Hence we have 
$$B_{j}^T{(B_{j}D_{j+1}B_{j}^T)}^{-1}B_{j}D_{j+1}
+C_{j}^T{(C_{j}D_{j+1}C_{j}^T)}^{-1}C_{j}D_{j+1}=I
$$
or
$$B_{j}^T{(B_{j}D_{j+1}B_{j}^T)}^{-1}B_{j}+C_{j}^T{(C_{j}D_{j+1}C_{j}^T)}^{-1}C_{j}
=D_{j+1}^{-1}.$$
which is (\ref{key}) and hence $h_j+ g_j=g_{j+1}$.

\section{Numerical Experiments}
We have implemented the prewavelet method for numerical solution of Poisson equations
over rectangular domains in MATLAB. We would like to demonstrate that our prewavelet
method is more efficient than the standard FEM method. 

In the following we provide three tables of
CPU times for numerical solutions based on 
our prewavelet method and the standard finite element method for various levels of 
refinement of an initial triangulation ($\Gamma_0$ which consists of two triangles) of 
the standard domain $[0, 1]\times [0, 1]$.

Let $V_j$ be the continuous linear finite element space over triangulation $\Gamma_j$
which is the $j$th refinement of $\Gamma_0$. For a test function $u$ which is the exact
solution of Poisson equation (\ref{PDE}), the finite element method is to compute
$u_j\in V_j$ directly while our prewavelet method computes $u_j$ by computing 
$w_k, k=1, \cdots, j$, i.e., $u_j=u_1+w_1+\cdots+ w_{j-1}$.  

In the following we present three tables of CPU times for computing numerical solutions
$u_j, j=4, 5, 6$ for three test solutions by using these two methods. Note that we
use the direct method coded in MATLAB to solve the associated linear equations. We 
shall present tables of CPU times based on Conjugate Gradient Method for the systems 
of equations next.

For an exact solution $u(x,y)=\sin(2\pi x)\sin(2\pi y)$ which clearly satisfies the
zero boundary conditions,  we list CPU times for computing numerical solutions 
$u_j, j=4, 5, 6$  by using these two methods in Table 1. 

\bigskip
\centerline{Table 1. CPU times to compute $u_j$ by the two methods} 
\begin{center}
\begin{tabular}{|c|c|c|}
\hline
    &FEM method&Prewavelet Method\\ 
\hline 
j=4&0.164531 seconds&0.204067 seconds\\ 
\hline
j=5&0.593587 seconds&0.519293 seconds\\ 
\hline
j=6&13.960323 seconds&6.222679 seconds\\ 
\hline
\end{tabular}
\end{center}

For an exact solution  $u(x,y)=xy(1-x)(1-y)$, the CPU times for numerical solutions by
these two methods are given in Table 2.

\bigskip
\centerline{Table 2. CPU times for computing $u_j$ by the two methods}
\begin{center}
\begin{tabular}{|c|c|c|}
\hline
CPU time&FEM method&Prewavelet Method\\ 
\hline 
j=4&0.150836 seconds&0.218282 seconds\\ 
\hline
j=5&0.574085 seconds&0.558071 seconds\\ 
\hline
j=6&13.896825 seconds&6.202557 seconds\\ 
\hline
\end{tabular}
\end{center}

We list the CPU times for computing numerical solutions $u_j, j=4, 5, 6$ of 
$u(x,y)=xy(1-x)(1-y)e^{8xy}$ by using these two methods in Table 3.

\bigskip
\centerline{Table 3. CPU times for computing $u_j$ by the two methods}
\begin{center}
\begin{tabular}{|c|c|c|}
\hline
CPU time&FEM method&Prewavelet Method\\ 
\hline 
j=4&0.144159 seconds&0.186389 seconds\\ 
\hline
j=5&0.584828 seconds&0.459181 seconds\\ 
\hline
j=6&13.877403 seconds&6.139101 seconds\\ 
\hline
\end{tabular}
\end{center}
It is clear from these three tables that the prewavelet method is much more efficient. 

Next we use the Conjugate Gradient Method to solve the linear systems associated with 
FEM. Let us consider iterative solution to $u_j$ for $ j=6$ with various accuracy.  
First let us consider the exact solution $u(x,y)=\sin(2\pi x)\sin(2\pi y)$.

\bigskip
\centerline{Table 4. CPU times for approximating the FEM solution $u_6$ by 
Conjugate Gradient Method} 
\begin{center}
\begin{tabular}{|c|c|c|}
\hline
$\epsilon$ &CPU times\\ 
\hline 
$10^{-8}$&5.411852 seconds\\ 
\hline
$10^{-9}$&5.783497 seconds\\ 
\hline
$10^{-10}$&6.221683 seconds\\ 
\hline
$10^{-11}$&6.616816 seconds\\ 
\hline
$10^{-12}$&6.917468 seconds\\ 
\hline
$10^{-13}$&7.836775 seconds\\ 
\hline
\end{tabular}
\end{center}

To approximate the FEM solution $u_6$ of the exact solution $u(x,y)=xy(1-x)(1-y)$ 
by the Conjugate Gradient Method, we list the CPU times in Table 5. 

\bigskip
\centerline{Table 5. CPU times for approximating the FEM solution $u_j$ by 
Conjugate Gradient Method} 
\begin{center}
\begin{tabular}{|c|c|c|}
\hline
$\epsilon$ &CPU times\\ 
\hline 
$10^{-8}$&4.476794 seconds\\ 
\hline
$10^{-9}$&4.878259 seconds\\ 
\hline
$10^{-10}$&5.306747 seconds\\ 
\hline
$10^{-11}$&5.887849 seconds\\ 
\hline
$10^{-12}$&6.811317 seconds\\ 
\hline
$10^{-13}$&6.754465 seconds\\ 
\hline
\end{tabular}
\end{center}
Finally let us consider the CPU times to approximate the FEM solution $u_6$ of 
$u(x,y)=xy(1-x)(1-y)e^{8xy}$ by the Conjugate Gradient Method.

\bigskip
\centerline{Table 6. CPU times for approximating the FEM solution by
 Conjugate Gradient Method} 
\begin{center}
\begin{tabular}{|c|c|c|}
\hline
$\epsilon$ &CPU times\\ 
\hline 
$10^{-8}$&10.110517 seconds\\ 
\hline
$10^{-9}$&10.740035 seconds\\ 
\hline
$10^{-10}$&11.319618  seconds\\ 
\hline
$10^{-11}$&11.810142 seconds\\ 
\hline
$10^{-12}$&12.320903 seconds\\ 
\hline
$10^{-13}$&13.103407 seconds\\ 
\hline
\end{tabular}
\end{center}
It is clear from all six tables, if we want an accurate iterative solution of $u_6$ 
within $10^{-12}$,  the prewavelet method appears better. 


\end{document}